\documentclass{article}
\usepackage{amsthm}
\usepackage{amsmath}
\usepackage{amsfonts}
\usepackage{amssymb}
\usepackage{calc}

\usepackage[all,ps,dvips,arc,knot]{xy}     % First choice
\def \R{{\mathbb R}}

\def \Z{{\mathbb Z}}
\def \A{{\mathcal A}}
\def \P{{\mathcal P}}
\def \F{{\mathcal F}}

\def \D{{\mathcal D}}

\def \OInt{\overrightarrow{\cap}}
\def \IntNum#1{<\! #1 \!>}
\def \bracket#1{\bigl< #1 \bigr>}
\def \Bracket#1{\Biggl<\ #1\ \Biggr>}
\def \bf{\bfseries}
\def \vecfrac#1#2{\genfrac{[}{]}{0pt}{}{#1}{#2}}

\def \ur#1{\,\overline{\,#1\vphantom{ly}\,}\raisebox{-.81ex}{\rule{.4pt}{2.8ex}}\,}
\def \lr#1{\,\underline{\,#1\vphantom{ly}\,}\raisebox{-.81ex}{\rule{.4pt}{2.8ex}}\,}
\def \ul#1{\,\raisebox{-.81ex}{\rule{.4pt}{2.8ex}}\overline{\ \vphantom{ly}#1\,}\,}
\def \ll#1{\,\raisebox{-.81ex}{\rule{.4pt}{2.8ex}}\underline{\ \vphantom{ly}#1\,}\,}
\def \Bur#1{\,\overline{\,#1\vphantom{\raisebox{-.90ex}{\rule{.4pt}{2.94ex}}}\,}\raisebox{-1.28ex}{\rule{.4pt}{3.72ex}}\,}
\def \Blr#1{\,\underline{\,#1\vphantom{\raisebox{-.90ex}{\rule{.4pt}{2.94ex}}}\,}\raisebox{-1.28ex}{\rule{.4pt}{3.72ex}}\,}
\def \Bul#1{\,\raisebox{-1.28ex}{\rule{.4pt}{3.72ex}}\overline{\ \vphantom{\raisebox{-.90ex}{\rule{.4pt}{2.94ex}}}#1\,}\,}

\def \Bbur#1{\,\overline{\,#1\vphantom{\raisebox{-1.37ex}{\rule{.4pt}{3.85ex}}}\,}\raisebox{-1.75ex}{\rule{.4pt}{4.6ex}}\,}

\theoremstyle{plain}
\newtheorem{thm}{Theorem}[section]

\newtheorem{lemma}[thm]{Lemma}

\newtheorem{conj}[thm]{Conjecture}
\theoremstyle{definition}
\newtheorem{defn}[thm]{Definition}
\theoremstyle{remark}

\begin{document}
% \squarebrackets        % This changes the references to [] brackets from ().

\date{}
\title{\Large\bf On Filamentations and Virtual Knots}

\author{David Hrencecin \\
Department of Mathematics, Statistics \\
and Computer Science (m/c 249)    \\
851 South Morgan Street   \\
University of Illinois at Chicago\\
Chicago, Illinois 60607-7045\\
$<$dhren@math.uic.edu$>$ \\
\\
and \\
\\
Louis H. Kauffman\\
Department of Mathematics, Statistics \\
and Computer Science (m/c 249)    \\
851 South Morgan Street   \\
University of Illinois at Chicago\\
Chicago, Illinois 60607-7045\\
$<$kauffman@uic.edu$>$}

\maketitle

\section{Introduction} \label{S:Summary}
Virtual knot theory is a recent generalization of knot theory.
One motivation for studying virtual knots comes from the methods
of describing knots through the use of chord diagrams.  In
Section~\ref{S:Definitions} we give the basic definitions for
this point of view.  In particular, we define oriented chord
diagrams and arrow diagrams.  The definition for a virtual knot
is also given, though we recommend~\cite{LouVirt} as an introduction
for the reader who is not already familiar with virtual knot theory.
Flat virtual diagrams and their equivalence classes are introduced in
Section~\ref{S:VKs}.

In Section~\ref{S:Filamentations}, we define a filamentation on a
Gauss chord diagram.  Filamentations were first introduced by
Scott Carter as a tool for detecting when an immersed curve can be
bounded by a disk.  We show that they can also detect when virtual knots
are non-trivial.  We demonstrate that a virtual knot diagram cannot be
reduced to a classical knot whenever a filamentation does not exist on its
chord diagram.  This result is Theorem~\ref{T:R_fil}, which we
prove from our combinatorial point of view.  In fact, we prove that
whenever there is no filamentation on a chord diagram, then any
associated flat virtual diagram cannot be reduced to a classical diagram.
There are many open problems in the classification
of flat virtual knots.

One interesting result of Theorem~\ref{T:R_fil} is that we have
found an infinite class of chord diagrams for which no
filamentation exists.  Section~\ref{S:Examples} explores this
infinite class of virtual knots ($K_n$) related to the chord
diagrams. The generalizations of both the Jones polynomial and the
fundamental group to virtual knots cannot detect any of these
$K_n$.  Using filamentations, we are able to show that each $K_n$
is non-trivial, although this method cannot distinguish between the $K_n$.
In Section~\ref{S:Biquandle}, we prove that the $K_n$ are an infinite
class of mutually distinct knots.

We would like to thank Scott Carter, Heather Dye and David Radford for
their very helpful conversations in the course of preparing this paper.

\section{Definitions} \label{S:Definitions}

We define a knot $K$ in the combinatorial sense, as a class of
diagrams which represent a generic projection of an embedding $S^1\to S^3$
or $S^1\sqcup\dots\sqcup S^1\to S^3$.  Each circle represented in such a diagram
is called a component, and if there is more than one component, the diagram (or related class)
is sometimes referred to as a link.
A strand in a knot diagram is the projection of a connected interval in the
embedded curve.
In a knot diagram,
Certain arcs of the projection of the curve embedded
in space are eliminated to form a knot diagram, creating broken strands
that indicate the over and under crossings.
At each crossing in the diagram there is an indicated over strand and a broken
under strand.
In this sense, there are two local strands at any crossing
in a diagram.  Any time we refer to a strand of a crossing, we mean one
of these two local strands.

\begin{defn}
Two knot diagrams $K$ and $K'$ are said to be \emph{equivalent} when there is a
finite sequence of the following (Reidemeister) moves which transform $K$ into $K'$:
% XXX -- beginning of pictures
$$
\hbox{\bf(R1)\ }
\xy 0;<1pc,0pc>:
0;(1.5,1.125)**\crv{(1.5,-0.75)&(2.25,1.075)}
,(1.5,1.125);(1.2,0.25)**\crv{(0.750,1.15)}
,(1.6,-0.1);(3,0)**\crv{(2.3,-0.50)}
\endxy \quad \xy 0;<1pc,0pc>:
\ar @{^{<}-_{>}} ,(0,0.25);(1,0.25)
\endxy \quad
\xy 0;<1pc,0pc>:*=dir{}
,0;(3,0)**\crv{(1.125,-0.75)&(1.875,0.75)}
\endxy \qquad\qquad
\hbox{\bf(R2)\ }
\xy 0;<1pc,0pc>:
,\ar @{-} @/_4pt/,(0,-1.5);(1.15,0)
,\ar @{-} @/_4pt/,(1.15,0);(0,1.5)
,\ar @{-} @/^4pt/|(0.45)\hole,(1.5,-1.5);(0.35,0)
,\ar @{-} @/^4pt/|(0.55)\hole,(0.35,0);(1.5,1.5)
\endxy \quad \xy 0;<1pc,0pc>:
\ar @{^{<}-_{>}} ,(0,0.25);(1,0.25)
\endxy \quad \xy 0;<1pc,0pc>:*=dir{}
\ar @{-} @/_/ ,(0,-1.5);(0,1.5)
\ar @{-} @/^/ ,(1.5,-1.5);(1.5,1.5)
\endxy
$$
$$
\hbox{\bf(R3)\ }
\xy 0;<1.5pc,0pc>:
a(0)="a1",a(60)="b1",a(120)="c1",a(180)="a2",a(240)="b2",a(300)="c2"
,\ar @{-}@/^1.5pt/|(0.69)\hole "a1";(0,-0.3)
,\ar @{-}@/^1.5pt/|(0.31)\hole (0,-0.3);"a2"
,\ar @{-}@/^3pt/ "b2";"b1"
,\ar @{-}@/^3pt/|(0.37)\hole "c1";"c2"
\endxy \quad \xy
0;<1pc,0pc>:\ar @{^{<}-_{>}} ,(0,0.25);(1,0.25)
\endxy \quad \xy 0;<1.5pc,0pc>:
a(0)="a1",a(60)="b1",a(120)="c1",a(180)="a2",a(240)="b2",a(300)="c2"
,\ar @{-}@/_1.5pt/|(0.65)\hole "a1";(0,0.3)
,\ar @{-}@/_1.5pt/|(0.35)\hole (0,0.3);"a2"
,\ar @{-}@/_4pt/ "b2";"b1"
,\ar @{-}@/_4pt/|(0.67)\hole "c1";"c2"
\endxy
$$
% XXX -- end of pictures
A \emph{knot} is an equivalence class of knot diagrams under the Reidemeister moves.
\end{defn}
A knot can have an orientation.  This is indicated on a knot diagram by drawing an arrowhead
on one or more strands in such a way that each component has a consistent labelling.
The Reidemeister moves for an oriented knot are the same as for unoriented knots.
We are free to apply the moves without paying attention to the particular orientations
of the strands involved.

There is also the notion of more than one kind of Reidemeister equivalence.
\emph{Regular isotopy} is equivalence under only the R2 and R3 moves.
\emph{Ambient isotopy} is equivalence under all three moves.
It is useful to distinguish between the two because there are invariants
which only cover regular isotopy.  For a more general treatise on knot theory,
see~\cite{LouOnKs} and~\cite{Kawauchi}.

In the study of Vassiliev (or finite type) knot invariants,
chord diagrams and weight systems have been used as a calculational
tool~\cite{Dror,GPV,PV}.
We will examine chord diagrams further and investigate their usefulness
in the general theory of knots.
Briefly, a chord diagram is a circle
(or set of disjoint circles\footnote{For links, we use a circle
to represent each component.  Each crossing involving a single
component is represented by a chord contained in the interior of that
component's circle.  All crossings between two different components
are represented by chords between each circles' exterior.})
with pairs of points on it, where each pair of points is connected by
a line segment, or \emph{chord}, in the interior of the circle.
We will assume the convention that the circle in a
chord diagram is oriented in a counterclockwise direction.
One way of thinking of a chord diagram is to view the outer circle(s) as the
pre-image of the projection of a knot.
\begin{defn}
The \emph{universe} (or shadow)~\cite{LouVirt} of a knot is a
generic projection of the embedding without specified over or under crossings.
We sometimes refer to this as a flattened knot diagram.
\end{defn}
The circle(s) in a chord diagram can be interpreted as the domain of a knot projection.
Each circle corresponds to a component in an associated embedding.
Hence we will sometimes refer to each circle as a component of the chord diagram.
Each chord connects the double points corresponding to a particular crossing.
This does not encode the types of crossings involved.  Additional structure is needed
in order to describe a knot diagram completely.

One solution to the problem of encoding the knot diagram
 is to add a sign and a direction to each chord~\cite{GPV,Polyak,PV}.
\begin{defn}
A \emph{signed arrow diagram} is a chord diagram in which each chord is
decorated with an arrow and given a sign.  The sign determines the crossing
orientation and the arrow points to the chord endpoint which lies on the
undercrossing strand in a related knot diagram.
\end{defn}
\begin{defn}
Dropping the arrow directions from a signed arrow diagram leaves a chord diagram
with a single sign on each chord. We call this a \emph{signed chord diagram}.
\end{defn}
There is a slight problem with developing a theory using signed arrow
diagrams.  We would like to have a chord diagram analogue to oriented knot universes
which generalizes signed arrow diagrams.  However in this case, the arrow and sign
information are mostly interdependent.
Signed chord diagrams do not store the same kind of information that
universes do.  It turns out that the Jones polynomial depends solely
on the underlying signed chord diagram of a knot, and this certainly
cannot be said for knot universes.  Each universe covers
multiple knot classes, many of which have different Jones polynomials.

%  If we drop the arrows, leaving a
%single sign on each chord, we do not have enough information to encode an
%oriented universe.  We will take care of this problem by defining an
%oriented chord diagram.

Our enhancements will modify signed arrow diagrams in a way that
allows us to generalize to oriented knot universes.
We start by defining an \emph{oriented chord diagram} and then give a
new definition of an \emph{arrow diagram} which builds on the underlying
structure of an oriented chord diagram.

\subsection{Oriented Chord Diagrams} \label{S:OCDs}
\begin{defn} \label{D:ocd}
An \emph{oriented chord diagram,} or OCD, is a chord diagram with a labelling
of~`$+$' or~`$-$' on each chord endpoint, so that each chord connects points
of opposite sign.
\end{defn}

\begin{figure}
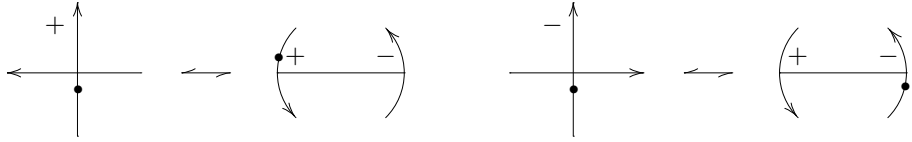

$$
\xy 0;/r2pc/:
{\ar (2.1,0);(0,0) \ar |<(0.35) *\dir{*} (1.1,-1);(1.1,1.1)}
,c+(-0.25,0.75) *{+}
,(5.25,0)="c" *\cir<2pc>{ur^ul} *\cir<2pc>{dl^dr}
,"c"+a(43)*{\object:a(41){\dir{>}}}, "c"-a(43)*{\object:a(220){\dir{>}}}
,{\ar@{-} "c"-(1,0)="a";"c"+(1,0)="b" \ar @{^{<}-_{>}} (2.75,0);(3.5,0) }
,"a"*+!DL{+},"a"+(0.02,0.25)*\dir{*} ,"b"*+!DR{-}
\endxy        \qquad\qquad
\xy 0;/r2pc/:
,{\ar (0,0);(2.1,0) \ar |<(0.35) *\dir{*} (1,-1);(1,1.1) }
,c+(-0.33,0.75) *{-}
,(5.25,0)="c" *\cir<2pc>{ur^ul} *\cir<2pc>{dl^dr}
,"c"+a(43)*{\object:a(41){\dir{>}}}, "c"-a(43)*{\object:a(220){\dir{>}}}
,{\ar@{-} "c"-(1,0)="a";"c"+(1,0)="b" \ar @{^{<}-_{>}} (2.75,0);(3.5,0) }
,"a"*+!DL{+},"b"-(0.02,0.22)*\dir{*},"b"*+!DR{-}
\endxy
$$
% XXX -- End picture
\caption{\label{Fi:OCDconv} Crossing conventions for oriented chord diagrams}
\end{figure}

An OCD encodes the universe of a knot diagram.  A neighborhood
of a chord endpoint (restricted to the circle) corresponds to one strand
of a crossing in a universe, and the sign on each chord endpoint determines
the local orientation relative to that strand.
In order to label a chord endpoint, start with its corresponding strand in the universe
and view that strand in the direction induced by the counter-clockwise orientation
of the circle.  Now look at the other strand of the crossing.
If it passes from right to left, label the
current chord endpoint with a~`$+$'.  Otherwise label it~`$-$'.
(See Figure~\ref{Fi:OCDconv}).
% XXX -- beginning of pictures
In other words, we view the two strands as vectors in the plane of projection.
The orientation comes from a standard right handed convention relating the
first vector to the second.  Switching the point of view from one chord endpoint
to the other reverses the local orientation.

\begin{figure}
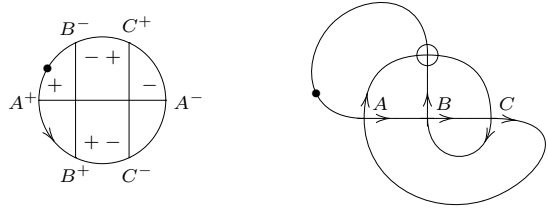

$$
\hbox{
\def\objectstyle{\scriptstyle}
\def\labelstyle{\scriptstyle}
\xy 0;<2pc,0pc>:
\POS(0,0)
*{\xybox{*\cir(1,0){} ,a(150)*{\bullet},a(-140)*\object:a(-144.5){\dir{>}}
,a(65)="a" *+!UR{+} *+!D{\phantom{-}C^+} ,a(115)="b" *+!UL{-} *+!D{B^-}
,a(180)="c" *+!DL{+} *!R{A^+} ,a(245)="d" *+!DL{+} *+!U{B^+}
,a(295)="e" *+!DR{-} *+!U{\phantom{-}C^-} ,(1,0)="f" *+!DR{-} *+!L{A^-}
, \ar @{-},"a";"e", \ar @{-},"d";"b", \ar @{-},"c";"f"
}\quad \ }="d1"
\POS(4.5,0)
*{\quad \ \xybox{0;<2pc,0pc>: (0,1.2)="s1", (0,4)="s2"
,(1,0)="A"*+!DL{A}, (2,0)="B"*+!DL{B}, (3,0)="C"*+!DL{C} ,(2,1)="V"*+\cir<4pt>{}
,(4.5,0)="a1",(3.0,-2.1)="a2",(1,-1)="a3" ,(1,1)="b1" ,(3,1)="c1"
,(3,-0.8)="d1",(2,-0.8)="d2" ,(2,2)="e1",(0,2)="e2",(0,0)="e3"
,"C";"A" **\crv{"a1"&"a2"&"a3"} ?(0.05)*\object:a(4){\dir{>}}
,"A";"V" **\crv{"b1"} ?(0.22)*\object:a(6.5){\dir{>}}
,"V";"C" **\crv{"c1"}
,"C";"B" **\crv{"d1"&"d2"} ?(0.19)*\object:a(6.5){\dir{>}}
,"V";"A" **\crv{"e1"&"e2"&"e3"} ?(0.8)*{\bullet}
,"A";"B" **\crv{} ?(0.4)*\dir{>} ,"B";"C" **\crv{} ?(0.4)*\dir{>}
,"B";"V" **\crv{} ?(0.4)*\dir{>}
}}="d2"
\endxy
}
$$
% XXX End picture
\caption{\label{Fi:OCDexample} An oriented chord diagram and related knot universe}
\end{figure}

Suppose we wish to recapture the universe from a single component OCD.
Choose a point on the circle that is not a chord endpoint.  Starting at this point,
follow around the circle in a counterclockwise direction.
While traversing the circle, draw a curve forming a diagram as follows:
\begin{itemize}
\item Each time the first endpoint of a new chord is encountered, draw a crossing
    obeying the sign convention and continue along the curve.
\item Whenever possible, draw the curve so that the only crossings which
    occur correspond to chords in the diagram.
\item When the second endpoint of a chord is encountered, connect the curve
    through the previously drawn crossing in the direction indicated by the sign.
\item Connect the endpoints of the curve when all the chords are accounted for.
\item If the curve is forced to cross itself at any point, circle the forced
    crossings to distinguish them from those associated with chords.
\end{itemize}
Figure~\ref{Fi:OCDexample} shows a knot universe along with the
OCD which represents it.
% XXX Begin picture
If there are multiple components in the OCD, proceed in the same fashion
for each one.

In the resulting universe, any crossing indicated by a chord in the OCD is
called a real (or classical) crossing.  The remaining crossings -- those forced
by the planar configuration of real crossings -- are called virtual.
Virtual crossings are like edge crossings in a drawing of a non-planar graph;
they do not exist in the original chord diagram, but are artifacts
of drawing the associated planar universe.
\begin{defn}
When we refer to crossing \emph{type}, we mean the distinction between real and virtual.
\end{defn}
In addition to diagrammatic representation, an OCD can be encoded in terms
of the chord endpoints and how they are ordered
while traversing the circle(s).  Suppose $\D$ is an OCD with $n$ chords,
$\{x_1, x_2,\dots, x_n\}$.  For each chord $x_i$, label the positive endpoint
$X_i^+$ and the
negative\footnote{We also use the negative superscript as an orientation
reversing operator: $(A^-)^- = A^+$.}
endpoint $X_i^-$.  In this paper, we will stick to the convention
of using lowercase letters when referring to chords and the corresponding
uppercase letters when referring to endpoints.  If $c$
is a chord in $\D$, then $C^+$ and $C^-$ are the positively and
negatively (resp.) oriented endpoints of~$c$.

To encode an oriented chord diagram $\D$, begin at any nonsingular point on
the circle.  Traverse the circle $\D$ in a counterclockwise direction and write down
the appropriate label for each endpoint as you pass it, until you return
to your original point on the circle.  Multiple components generate
multiple codes: one for each component.
Consider the OCD in Figure~\ref{Fi:OCDexample}.  An example of a code associated with it is
$A^+B^+C^{-}A^{-}C^+B^{-}$.
These codes are unique up to cyclic permutation
(to account for where you start on the circle),
along with any permutation of the chord labels.

That is,
\[
A^+B^+C^{-}A^{-}C^+B^{-},\quad
C^{-}A^{-}C^+B^{-}A^+B^+,\quad \hbox{and}\quad
A^+C^+B^{-}A^{-}B^+C^{-}
\]
are all codes for the same OCD.

\subsection{Arrow Diagrams} \label{S:ADs}
So far, we have described OCD's and their relationship with knot universes.
We define an arrow diagram within this framework.
%If the chord diagram comes from a classical knot diagram's universe, we can
%always go between the universe and an OCD.

\begin{defn}
An \emph{arrow diagram,} or AD, is an oriented chord diagram in which
the chords are replaced with arrows.  The direction of each arrow gives the
extra structure on a crossing in the related knot universe by the convention
that the arrow points from the overcrossing to the undercrossing strand.
\end{defn}

\begin{figure}
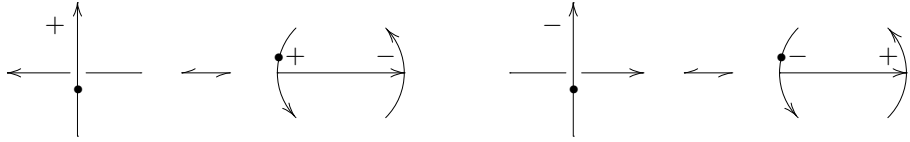

% XXX -- beginning of pictures
$$
\xy 0;/r2pc/:
{\ar |<(0.47)\hole (2.1,0);(0,0) \ar |<(0.35) *\dir{*} (1.1,-1);(1.1,1.1) }
,c+(-0.25,0.75) *{+}
,(5.25,0)="c" *\cir<2pc>{ur^ul} *\cir<2pc>{dl^dr}
,"c"+a(43)*{\object:a(41){\dir{>}}},"c"-a(43)*{\object:a(220){\dir{>}}}
,{\ar "c"-(1,0)="a";"c"+(1,0)="b" \ar @{^{<}-_{>}} (2.75,0);(3.5,0) }
,"a"*+!DL{+},"a"+(0.02,0.25)*\dir{*} ,"b"*+!DR{-}
\endxy        \qquad\qquad
\xy 0;/r2pc/:
,{\ar |<(0.48)\hole (0,0);(2.1,0) \ar |<(0.35) *\dir{*} (1,-1);(1,1.1) }
,c+(-0.33,0.75) *{-}
,(5.25,0)="c" *\cir<2pc>{ur^ul} *\cir<2pc>{dl^dr}
,"c"+a(43)*{\object:a(41){\dir{>}}},"c"-a(43)*{\object:a(220){\dir{>}}}
,{\ar "c"-(1,0)="a";"c"+(1,0)="b" \ar @{^{<}-_{>}} (2.75,0);(3.5,0) }
,"a"*+!DL{-},"a"+(0.02,0.25)*\dir{*},"b"*+!DR{+}
\endxy
$$
% XXX -- End picture
\caption{\label{Fi:ADconv} Crossing conventions for arrow diagrams}
\end{figure}
Note that the local orientation on the base endpoint of an arrow
(corresponding to the overcrossing strand) matches the classical
knot theoretical convention of crossing sign for an oriented knot
(see Figure~\ref{Fi:ADconv}).  The local orientation of one chord endpoint
is always opposite to the other endpoint, so we will only label
the base endpoint of each arrow in an AD.  On the surface, an arrow diagram
looks just like a signed arrow diagram.
The difference is that the arrow diagram sign refers specifically to the
local orientation of the arrow basepoint, rather than the sign
of the arrow itself.

For the purposes of this paper, the code associated with an arrow diagram
will be the same as the code associated with the underlying OCD.  If more
information is needed, we can use the subscripts `$o$' and `$u$' to denote
arrow basepoints and endpoints (respectively).  For example,
$A^+_oB^+_uC^{-}_uA^{-}_uC^+_oB^{-}_o$ encodes an arrow diagram with the
same underlying code as the OCD in Figure~\ref{Fi:OCDexample}.
We emphasize again that this is not the
same as the convention for signed arrow diagrams where
only one sign is associated with each arrow.

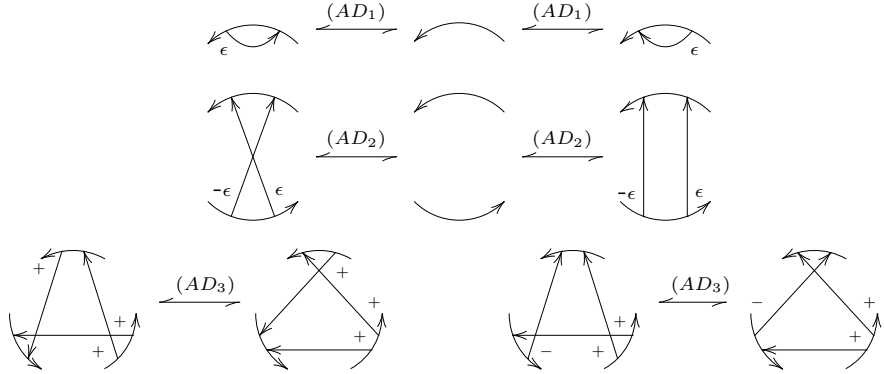
\begin{figure}
\[
{ \def\objectstyle{\scriptstyle}
\def\labelstyle{\scriptstyle}
\begin{gathered}
\xymatrix@C=30pt{
*!(0,5){\xybox{;<24pt,0pt>:0 *\cir(1,0){ul^dl} ,a(135)*{\object:a(132){\dir{>}}}
    ,(.383,.924)="a", (-.383,.924)="b"
    ,\ar_<(0.15){\epsilon} @/_/ ,a(115);a(65)}}
    \ar@{^{<}-_{>}}[r]^{(AD_1)}
&*!(0,5){\xybox{;<24pt,0pt>:0 *\cir(1,0){ul^dl} ,a(135)*{\object:a(222){\dir{>}}} }}
    \ar@{^{<}-_{>}}[r]^{(AD_1)}
&*!(0,5){\xybox{;<24pt,0pt>:0 *\cir(1,0){ul^dl} ,a(135)*{\object:a(222){\dir{>}}}
    ,(.383,.924)="a", (-.383,.924)="b"
    ,\ar^<(0.15){\epsilon} @/^/ ,a(65);a(115)}}
}\\
\xymatrix@C=30pt{
*{\xybox{;<24pt,0pt>:0 *\cir(1,0){ul^dl} *\cir(1,0){dr^ur}
    ,a(135)*{\object:a(132){\dir{>}}} ,a(315)*{\object:a(312){\dir{>}}}
    ,\ar^<(0.18){{\hbox{-}}\epsilon} ,a(250);a(70), \ar_<(0.18){\epsilon} ,a(290);a(110)}}
    \ar@{^{<}-_{>}}[r]^{(AD_2)}
&*{\xybox{;<24pt,0pt>:0 *\cir(1,0){ul^dl} *\cir(1,0){dr^ur}
    ,a(135)*{\object:a(222){\dir{>}}} ,a(315)*{\object:a(42){\dir{>}}} }}
\ar@{^{<}-_{>}}[r]^{(AD_2)}
&*{\xybox{;<24pt,0pt>:0 *\cir(1,0){ul^dl} *\cir(1,0){dr^ur}
    ,a(135)*{\object:a(211){\dir{>}}} ,a(315)*{\object:a(42){\dir{>}}}
    ,\ar^<(0.2){{\hbox{-}}\epsilon} ,a(250);a(110), \ar_<(0.2){\epsilon} ,a(290);a(70)}}
}\\
\xymatrix@C=30pt{
*{\xybox{;<24pt,0pt>:0
    ,a(300)="a";"a"+a(30),**{},a(0) *{\object:a(57){\dir{>}}},{\ellipse^{}}
    ,a(180)="c";"c"+a(270),**{},a(240)="d" *{\object:a(57){\dir{>}}},{\ellipse^{}}
    ,a(60)="e";"e"+a(150),**{},a(120)="f" *{\object:a(57){\dir{>}}},{\ellipse^{}}
    ,{\ar_<(0.2){_{+}} a(100);a(225),\ar_<(0.1){_{+}} a(340);a(200),\ar^<(0.1){_{+}} a(315);a(80) }
    }}
\ar@{^{<}-_{>}}[r]^{(AD_3)}
&\ {\xybox{;<24pt,0pt>:0
    ,a(300)="a";"a"+a(30),**{},a(0) *{\object:a(57){\dir{>}}},{\ellipse^{}}
    ,a(180)="c";"c"+a(270),**{},a(240)="d" *{\object:a(57){\dir{>}}},{\ellipse^{}}
    ,a(60)="e";"e"+a(150),**{},a(120)="f" *{\object:a(57){\dir{>}}},{\ellipse^{}}
    ,{\ar^<(0.1){_{+}} a(75);a(200),\ar_<(0.1){_{+}} a(325);a(215),\ar_<(0.25){_{+}} a(340);a(105) }
    }}
    }
\hspace{0.5in} \xymatrix@C=25pt{
*{\xybox{;<24pt,0pt>:0
    ,a(300)="a";"a"+a(30),**{},a(0) *{\object:a(57){\dir{>}}},{\ellipse^{}}
    ,a(180)="c";"c"+a(270),**{},a(240)="d" *{\object:a(57){\dir{>}}},{\ellipse^{}}
    ,a(60)="e";"e"+a(150),**{},a(120)="f" *{\object:a(57){\dir{>}}},{\ellipse^{}}
    ,{\ar_<(0.1){_{-}} a(225);a(100),\ar_<(0.1){_{+}} a(340);a(200),\ar^<(0.1){_{+}} a(315);a(80) }
    }}
\ar@{^{<}-_{>}}[r]^{(AD_3)}
&\ {\xybox{;<24pt,0pt>:0
    ,a(300)="a";"a"+a(30),**{},a(0) *{\object:a(57){\dir{>}}},{\ellipse^{}}
    ,a(180)="c";"c"+a(270),**{},a(240)="d" *{\object:a(57){\dir{>}}},{\ellipse^{}}
    ,a(60)="e";"e"+a(150),**{},a(120)="f" *{\object:a(57){\dir{>}}},{\ellipse^{}}
    ,{\ar^<(0.25){_{-}} a(200);a(75),\ar_<(0.1){_{+}} a(325);a(215),\ar_<(0.25){_{+}} a(340);a(105) }
    }}
}
\end{gathered}}
% XXX -- End picture
\]
\caption{\label{Fi:ADmoves} The arrow diagram moves}
\end{figure}
We can take the set of all arrow diagrams and define an
equivalence relation under the Reidemeister moves translated into
AD form (See~\cite{GPV}). Figure~\ref{Fi:ADmoves} lists these moves.
The $\epsilon$ on the arrows is meant to be either $+$ or $-$.
Each arc represents a portion of a circle in an arrow diagram, but
it is not necessary to assume that the relative
arc placements must be as shown.  In the $AD_2$ move, the two
arcs might lie on separate circles of a multiple component AD.  Further, in
the $AD_3$ move, we allow flexibility in the ordering of the arcs, even
within a single component diagram.
Simply put, the arc placements in a single application of an AD move can differ
from Figure~\ref{Fi:ADmoves}, provided that arc ordering does not change
across the move.

Consider the lower left version of the $AD_3$ move.  One application of
$AD_3$ might change an associated code's subcode from
\begin{align*}
A^-_uB^-_u\dots C^+_oA^+_o\dots C^-_uB^+_o \qquad &\hbox{to}\qquad
B^-_uA^-_u\dots A^+_oC^+_o\dots B^+_oC^-_u \\
\intertext{under the move. Another example of $AD_3$ might change the subcode from}
C^+_oA^+_o\dots A^-_uB^-_u\dots C^-_uB^+_o \qquad &\hbox{to}\qquad
A^+_oC^+_o\dots B^-_uA^-_u\dots B^+_oC^-_u
\end{align*}
in an arrow diagram.
With this in mind, these moves completely describe the reformulation of
Reidemeister moves into diagrammatic form and translate the
combinatorics of knot equivalence into arrow equivalence.

\section{Virtual Knots} \label{S:VKs}
It turns out that this new class of objects generalizes the classical knots.
Given an abstract arrow diagram $\A$, we find that
it is not always possible to realize $\A$ as a knot projection on a sphere
or plane.  Such a non-planar AD is a virtual knot.  A more diagrammatic
definition follows:
\begin{defn} \label{D:vk}
A \emph{virtual knot diagram} is a generic immersion $S^1\sqcup\dots\sqcup S^1\to \R^2$
such that each double point is labelled with either a real (over or under) or a virtual
(circled) crossing.  Thus a \emph{virtual knot} is a class of equivalent
virtual diagrams, where two virtual knot diagrams are said to be equivalent
when one can be transformed into the other by a finite sequence of real
Reidemeister moves (R1), (R2) and (R3), along with the following virtual moves:
% XXX -- beginning of pictures
$$
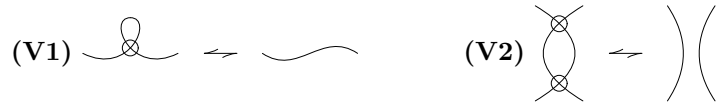
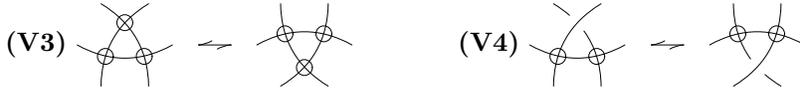
\hbox{\bf(V1)\ }
\xy 0;<1pc,0pc>:
0;(1.5,1.125)**\crv{(1.5,-0.75)&(2.25,1.125)} ?(0.45)*\cir<3pt>{}
,(1.5,1.125);(3,0)**\crv{(0.75,1.125)&(1.5,-0.75)}
\endxy \quad \xy 0;<1pc,0pc>:
\ar @{^{<}-_{>}} ,(0,0.25);(1,0.25)
\endxy \quad \xy
0;<1pc,0pc>:*=dir{}
,0;(3,0)**\crv{(1.125,-0.75)&(1.875,0.75)}
\endxy \qquad\qquad
\hbox{\bf(V2)\ }
\xy 0;<1pc,0pc>:
(0.75,0.94)*\cir<3pt>{}
,(0.75,-0.94)*\cir<3pt>{}
,\ar @{-} @/_15pt/,(0,-1.5);(0,1.5)
,\ar @{-} @/^15pt/,(1.5,-1.5);(1.5,1.5)
\endxy \quad \xy
0;<1pc,0pc>:
\ar @{^{<}-_{>}} ,(0,0.25);(1,0.25)
\endxy \quad \xy
0;<1pc,0pc>:*=dir{}
\ar @{-} @/_/ ,(0,-1.5);(0,1.5)
\ar @{-} @/^/ ,(1.5,-1.5);(1.5,1.5)
\endxy
$$
$$
\hbox{\bf(V3)\ }
\xy 0;<8.5pt,0pc>:
a(-30)*\cir<3pt>{},a(90)*\cir<3pt>{},a(210)*\cir<3pt>{}
,0;<1.5pc,0pc>:
a(0)="a1",a(60)="b1",a(120)="c1",a(180)="a2",a(240)="b2",a(300)="c2"
,\ar @{-}@/^5pt/,"a1";"a2"
,\ar @{-}@/^5pt/,"b2";"b1"
,\ar @{-}@/^5pt/,"c1";"c2"
\endxy \quad \xy
0;<1pc,0pc>:\ar @{^{<}-_{>}} ,(0,0.25);(1,0.25)
\endxy \quad \xy
0;<8.5pt,0pc>:
a(30)*\cir<3pt>{},a(150)*\cir<3pt>{},a(270)*\cir<3pt>{}
,0;<1.5pc,0pc>:
a(0)="a1",a(60)="b1",a(120)="c1",a(180)="a2",a(240)="b2",a(300)="c2"
,\ar @{-}@/_5pt/,"a1";"a2"
,\ar @{-}@/_5pt/,"b2";"b1"
,\ar @{-}@/_5pt/,"c1";"c2"
\endxy \qquad\qquad
\hbox{\bf(V4)\ }
\xy 0;<8.5pt,0pc>:
a(-30)*\cir<3pt>{},a(210)*\cir<3pt>{}
,0;<1.5pc,0pc>:
a(0)="a1",a(60)="b1",a(120)="c1",a(180)="a2",a(240)="b2",a(300)="c2"
,\ar @{-}@/^5pt/ "a1";"a2"
,\ar @{-}@/^5pt/ "b2";"b1"
,\ar @{-}@/^5pt/|(0.30)\hole "c1";"c2"
\endxy \quad \xy
0;<1pc,0pc>:\ar @{^{<}-_{>}} ,(0,0.25);(1,0.25)
\endxy \quad \xy
0;<8.5pt,0pc>:
a(30)*\cir<3pt>{},a(150)*\cir<3pt>{}
,0;<1.5pc,0pc>:
a(0)="a1",a(60)="b1",a(120)="c1",a(180)="a2",a(240)="b2",a(300)="c2"
,\ar @{-}@/_5pt/ "a1";"a2"
,\ar @{-}@/_5pt/ "b2";"b1"
,\ar @{-}@/_5pt/|(0.70)\hole "c1";"c2"
\endxy
$$
\end{defn}
% XXX -- end of pictures
As in the classical case, a virtual knot with multiple components is sometimes called
a virtual link.  For this paper, we will use the term `virtual knot' to refer to both
virtual knots and links.  When we do not wish to consider virtual crossings, we
will use the term `classical knot' to mean knots and links without virtual crossings.

An alternate definition to virtual Reidemeister equivalence is to allow the
classical Reidemeister moves with the addition of a more general ``detour move"
(See~\cite{LouVirt,LouVirt2}):
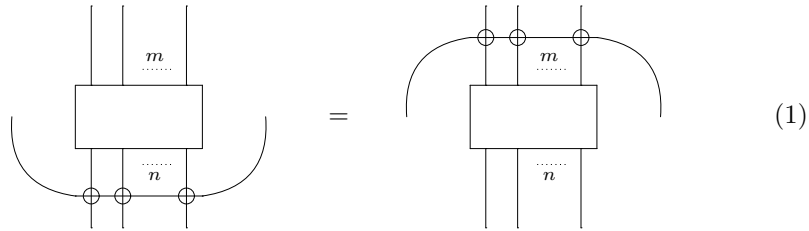
\begin{equation} \label{E:DetourMove}
\xy 0;/r1pc/:@={(0,1),(4,1),(4,-1),(0,-1)},s0="prev" @@{;"prev";**\dir{-}="prev"}
@i @={(.5,-2.5),(1.5,-2.5),(3.5,-2.5)} @@{*\cir<3pt>{}}
,\ar^{m}@{.} (2,1.5);(3,1.5),\ar_{n}@{.} (2,-1.5);(3,-1.5)
,\ar@{-} (.5,1);(.5,3.5),\ar@{-} (1.5,1);(1.5,3.5),\ar@{-} (3.5,1);(3.5,3.5)
,\ar@{-} (.5,-1);(.5,-3.5),\ar@{-} (1.5,-1);(1.5,-3.5),\ar@{-} (3.5,-1);(3.5,-3.5)
,\ar@{-} (0,-2.5);(4,-2.5),\ar@{-}@/_9pt/ (-2,0);(0,-2.5),\ar@{-}@/_9pt/ (4,-2.5);(6,0)
\endxy
\qquad = \qquad\quad
\xy 0;/r1pc/:@={(0,1),(4,1),(4,-1),(0,-1)},s0="prev" @@{;"prev";**\dir{-}="prev"}
@i @={(.5,2.5),(1.5,2.5),(3.5,2.5)} @@{*\cir<3pt>{}}
,\ar^{m}@{.} (2,1.5);(3,1.5),\ar_{n}@{.} (2,-1.5);(3,-1.5)
,\ar@{-} (.5,1);(.5,3.5),\ar@{-} (1.5,1);(1.5,3.5),\ar@{-} (3.5,1);(3.5,3.5)
,\ar@{-} (0,2.5);(4,2.5),\ar@{-}@/^9pt/ (-2,0);(0,2.5),\ar@{-}@/^9pt/ (4,2.5);(6,0)
,\ar@{-} (.5,-1);(.5,-3.5),\ar@{-} (1.5,-1);(1.5,-3.5),\ar@{-} (3.5,-1);(3.5,-3.5)
\endxy
\end{equation}
In the detour move, any number of strands may emanate from the top and bottom of the
tangle (represented by a box).  The idea is that in a virtual diagram, if we have an
arc with any number of consecutive virtual crossings, then we can cut that arc out
and replace it with another arc connecting the same points, provided that any crossings
on the new arc are also virtual.  It is easily seen that this yields the same equivalence.

\begin{defn}
If $K$ is a knot diagram, then $AD(K)$ is the arrow diagram \emph{related to} $K$.
\end{defn}
Notice that there are some mixed moves which are not allowed.  Consider
the following moves:
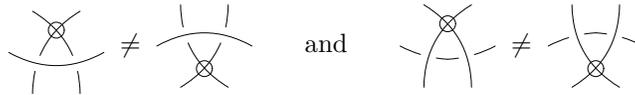
\[
\xy ,0;<1.5pc,0pc>:
<0pt,8.5pt>*\cir<3pt>{}
,a(0)="a1",a(60)="b1",a(120)="c1",a(180)="a2",a(240)="b2",a(300)="c2"
,\ar @{-}@/^5pt/"a1";"a2"
,\ar @{-}@/^5pt/|(0.30)\hole "b2";"b1"
,\ar @{-}@/^5pt/|(0.70)\hole "c1";"c2"
\endxy
\ \ne\
\xy ,0;<1.5pc,0pc>:
<0pt,-8.5pt>*\cir<3pt>{}
,a(0)="a1",a(60)="b1",a(120)="c1",a(180)="a2",a(240)="b2",a(300)="c2"
,\ar @{-}@/_5pt/"a1";"a2"
,\ar @{-}@/_5pt/|(0.70)\hole  "b2";"b1"
,\ar @{-}@/_5pt/|(0.30)\hole "c1";"c2"
\endxy
\qquad\hbox{and}\qquad
\xy ,0;<1.5pc,0pc>:
a(0)="a1",a(60)="b1",a(120)="c1",a(180)="a2",a(240)="b2",a(300)="c2"
,<0pt,8.5pt> *\cir<3pt>{},<0pt,-5pt>="a3"
,\ar @{-}@/^1pt/|(0.62)\hole "a1";"a3"
,\ar @{-}@/^1pt/|(0.39)\hole "a3";"a2"
,\ar @{-}@/^5pt/ "b2";"b1"
,\ar @{-}@/^5pt/ "c1";"c2"
\endxy
\ \ne\
\xy ,0;<1.5pc,0pc>:
a(0)="a1",a(60)="b1",a(120)="c1",a(180)="a2",a(240)="b2",a(300)="c2"
,<0pt,-8.5pt> *\cir<3pt>{},<0pt,5pt>="a3"
,\ar @{-}@/_1pt/|(0.62)\hole "a1";"a3"
,\ar @{-}@/_1pt/|(0.39)\hole "a3";"a2"
,\ar @{-}@/_5pt/ "b2";"b1"
,\ar @{-}@/_5pt/ "c1";"c2"
\endxy
\]
These are forbidden as knot diagrammatical analogues to
the AD moves.  They change the related arrow diagram in a manner which is
not equivalent under the AD moves; they permute two adjacent arrow
endpoints.  However, inclusion of the above left (overstrand) version
of this move has been studied in the form of welded braids~\cite{Fenn,FennEtal},
a generalization of braid theory preceding virtual knot theory.

From the definition, we note that a virtual knot or link can also
be classical.  This happens when it can be represented by a diagram
in which all of the crossings are real.  Further, since the virtual
moves (V1)-(V4) leave the related arrow diagram unchanged, they also
preserve the classical knot type.  However, it is possible to apply an
arrow diagram move to a chord diagram for a classical knot in such a way
that the resulting knot diagram is no longer classical.

There has been progress in applying well known knot invariants
to the virtual theory.  In~\cite{LouVirt} Kauffman extended the
fundamental group, the Jones polynomial and classes of quantum link invariants
to virtual knots and gave examples of non-trivial virtual knots with trivial
Jones Polynomial and trivial\footnote{by trivial fundamental
group, we mean a group isomorphic to the integers} fundamental group.
We will cover the fundamental group and the Jones polynomial
in Section~\ref{S:Examples}.

One of the results we will assume is the following, proved by
Kauffman along with Goussarov, Polyak and Viro:
\begin{lemma}
If $K$ and $K'$ are classical knots which are equivalent under virtual Reidemeister
equivalence, then they are also equivalent under classical Reidemeister equivalence.
\end{lemma}

We will also use the following result:
\begin{thm} \label{T:AD_Equiv}
Virtual equivalence and AD equivalence are the same thing.  That is,
if two arrow diagrams are equivalent under AD moves, then any virtual knots
related to those two diagrams will be virtually equivalent.  Likewise,
if $K\simeq K'$, then $AD(K)\simeq AD(K')$.
\end{thm}
\begin{proof}
This is a direct result of the definitions.
\end{proof}
An immediate question that arises is how to determine when a virtual knot
is classical.  In~\cite{LouVirt}, Kauffman introduces flat Reidemeister
equivalence as one approach to this problem.
\begin{defn}
A \emph{flat Reidemeister move} is a classical or virtual Reidemeister
move in which the over/under information at each real crossing is
suppressed, so that all that is important is the distinction between
crossing types (real and virtual).  Two universes
are \emph{flat (Reidemeister) equivalent} if there is a finite
sequence of flat Reidemeister moves taking the one universe to
the other.
An equivalence class of knot universes under flat
equivalence is called a \emph{flat knot.}
\end{defn}

\begin{figure}
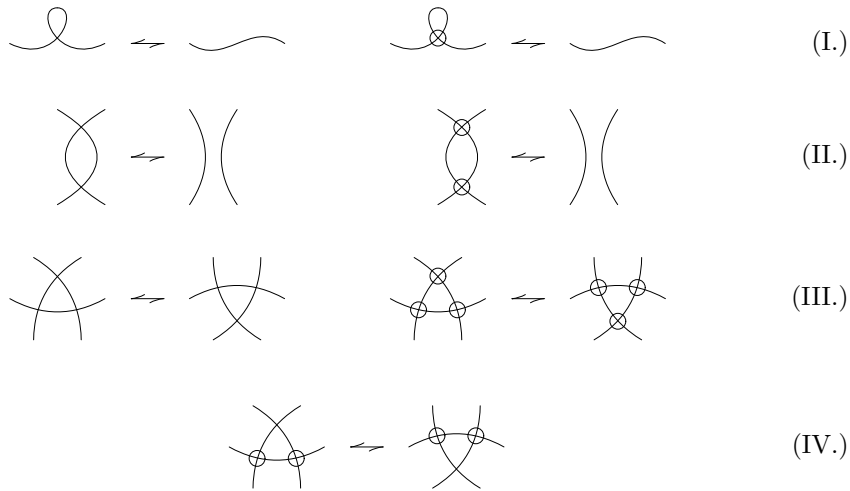

% XXX -- Begin picture
\begin{alignat*}{2}
{\xy
0;<1pc,0pc>:
0;(1.5,1.125)**\crv{(1.5,-0.75)&(2.25,1.125)}
,(1.5,1.125);(3,0)**\crv{(0.75,1.125)&(1.5,-0.75)}
\endxy} &\quad {\xy
0;<1pc,0pc>:\ar @{^{<}-_{>}} ,(0,0.25);(1,0.25)
\endxy} \quad {\xy
0;<1pc,0pc>:*=dir{}
,0;(3,0)**\crv{(1.125,-0.75)&(1.875,0.75)}
\endxy} &\qquad\qquad {\xy
0;<1pc,0pc>:
0;(1.5,1.125)**\crv{(1.5,-0.75)&(2.25,1.125)} ?(0.45)*\cir<3pt>{}
,(1.5,1.125);(3,0)**\crv{(0.75,1.125)&(1.5,-0.75)}
\endxy} &\quad {\xy
0;<1pc,0pc>:\ar @{^{<}-_{>}} ,(0,0.25);(1,0.25)
\endxy} \quad {\xy
0;<1pc,0pc>:*=dir{}
,0;(3,0)**\crv{(1.125,-0.75)&(1.875,0.75)}
\endxy}\qquad\tag{I.}
\\&&&\notag
\\
{\xy
0;<1pc,0pc>:
(0.75,0.94)
,(0.75,-0.94)
,\ar @{-} @/_15pt/,(0,-1.5);(0,1.5)
,\ar @{-} @/^15pt/,(1.5,-1.5);(1.5,1.5)
\endxy} &\quad {\xy
0;<1pc,0pc>:\ar @{^{<}-_{>}} ,(0,0.25);(1,0.25)
\endxy} \quad {\xy
0;<1pc,0pc>:*=dir{}
\ar @{-} @/_/ ,(0,-1.5);(0,1.5)
\ar @{-} @/^/ ,(1.5,-1.5);(1.5,1.5)
\endxy} &\qquad\qquad {\xy
0;<1pc,0pc>:
(0.75,0.94)*\cir<3pt>{}
,(0.75,-0.94)*\cir<3pt>{}
,\ar @{-} @/_15pt/,(0,-1.5);(0,1.5)
,\ar @{-} @/^15pt/,(1.5,-1.5);(1.5,1.5)
\endxy} &\quad {\xy
0;<1pc,0pc>:\ar @{^{<}-_{>}} ,(0,0.25);(1,0.25)
\endxy} \quad {\xy
0;<1pc,0pc>:*=dir{}
\ar @{-} @/_/ ,(0,-1.5);(0,1.5)
\ar @{-} @/^/ ,(1.5,-1.5);(1.5,1.5)
\endxy}\qquad\tag{II.}
\\&&&\notag
\\
{\xy 0;<8.5pt,0pc>:
,0;<1.5pc,0pc>:
a(0)="a1",a(60)="b1",a(120)="c1",a(180)="a2",a(240)="b2",a(300)="c2"
,\ar @{-}@/^5pt/,"a1";"a2"
,\ar @{-}@/^5pt/,"b2";"b1"
,\ar @{-}@/^5pt/,"c1";"c2"
\endxy} &\quad {\xy
0;<1pc,0pc>:\ar @{^{<}-_{>}} ,(0,0.25);(1,0.25)
\endxy} \quad {\xy
0;<8.5pt,0pc>:
,0;<1.5pc,0pc>:
a(0)="a1",a(60)="b1",a(120)="c1",a(180)="a2",a(240)="b2",a(300)="c2"
,\ar @{-}@/_5pt/,"a1";"a2"
,\ar @{-}@/_5pt/,"b2";"b1"
,\ar @{-}@/_5pt/,"c1";"c2"
\endxy} &\qquad\qquad {\xy
0;<8.5pt,0pc>:
a(-30)*\cir<3pt>{},a(90)*\cir<3pt>{},a(210)*\cir<3pt>{}
,0;<1.5pc,0pc>:
a(0)="a1",a(60)="b1",a(120)="c1",a(180)="a2",a(240)="b2",a(300)="c2"
,\ar @{-}@/^5pt/,"a1";"a2"
,\ar @{-}@/^5pt/,"b2";"b1"
,\ar @{-}@/^5pt/,"c1";"c2"
\endxy} &\quad {\xy
0;<1pc,0pc>:\ar @{^{<}-_{>}} ,(0,0.25);(1,0.25)
\endxy} \quad {\xy
0;<8.5pt,0pc>:
a(30)*\cir<3pt>{},a(150)*\cir<3pt>{},a(270)*\cir<3pt>{}
,0;<1.5pc,0pc>:
a(0)="a1",a(60)="b1",a(120)="c1",a(180)="a2",a(240)="b2",a(300)="c2"
,\ar @{-}@/_5pt/,"a1";"a2"
,\ar @{-}@/_5pt/,"b2";"b1"
,\ar @{-}@/_5pt/,"c1";"c2"
\endxy}\qquad\tag{III.}
\end{alignat*}
\smallskip
\begin{equation}
{\xy 0;<8.5pt,0pc>:
a(-30)*\cir<3pt>{},a(210)*\cir<3pt>{}
,0;<1.5pc,0pc>:
a(0)="a1",a(60)="b1",a(120)="c1",a(180)="a2",a(240)="b2",a(300)="c2"
,\ar @{-}@/^5pt/ "a1";"a2"
,\ar @{-}@/^5pt/ "b2";"b1"
,\ar @{-}@/^5pt/ "c1";"c2"
\endxy} \quad {\xy
0;<1pc,0pc>:\ar @{^{<}-_{>}} ,(0,0.25);(1,0.25)
\endxy} \quad {\xy
0;<8.5pt,0pc>:
a(30)*\cir<3pt>{},a(150)*\cir<3pt>{}
,0;<1.5pc,0pc>:
a(0)="a1",a(60)="b1",a(120)="c1",a(180)="a2",a(240)="b2",a(300)="c2"
,\ar @{-}@/_5pt/ "a1";"a2"
,\ar @{-}@/_5pt/ "b2";"b1"
,\ar @{-}@/_5pt/ "c1";"c2"
\endxy}\qquad\tag{IV.}
\end{equation}
% XXX -- End picture
\caption{\label{Fi:FlatMoves} The Flat Reidemeister Moves}
\end{figure}

In Figure~\ref{Fi:FlatMoves}, we illustrate the (flat) Reidemeister moves for flat
virtual diagrams.  Note that a sequence of virtual crossings can be detoured
across flat crossings (using the mixed move given by IV), but not vice versa.
The virtual detour move given in~\eqref{E:DetourMove} applies for flat knots
as well, and can be thought of as equivalent to move IV in the presence of
the other moves.

The remarkable fact about flat virtual diagrams is that while they are often non-trivial
(and hence non-classical), we have very few invariants at the present time which detect
and classify them.  For example, consider the following flat diagram:

% XXX -- Begin picture
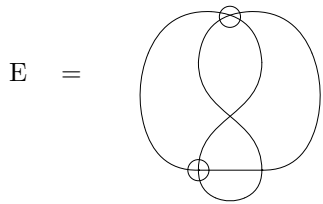
\[
\hbox{E}\quad
=\!\!\!\!
{\xy 0;/r2pc/:
(0.2,1)="a1",(0.8,1)="b1",(0.5,0.92)*{\cir<4pt>{}}
,(0,0.2)="a2",(1,0.2)="b2"
,(0,-1.5)="a3",(1,-1.5)="b3"
,(0,-1.5)="a4" *{\cir<4pt>{}},(1,-1.5)="b4"
,(0,-1.5)="a5",(1,-1.5)="b5",
,\ar@/^6pt/@{-} "a1";"b2",\ar@/_6pt/@{-} "b1";"a2"
,{\save "b2"="temp","a2";"a3":"temp"::(0,0);(1,1) **\crv{(0.5,0)&(0.5,1)},(0,1);(1,0) **\crv{(0.5,1)&(0.5,0)},(0.5,0.5) \restore}
,\ar@{-} "a4";"b4" ,\ar@{-} "a3";"a5" ,\ar@{-} "b3";"b5"
\hcap~{"b1"}{"b1"+(2.2,0.1)}{"b4"}{"b4"+(2,-.05)}
\hcap~{"a1"}{"a1"-(2.2,-0.1)}{"a4"}{"a4"-(2,0.05)}
\vcap~{"a5"-(0,1)}{"b5"-(0,1)}{"a5"}{"b5"}
\endxy}
\]
% XXX -- End picture
At the time of this writing, we do not have a proof that the flat diagram E given above
is inequivalent to a circle.

A simple example of a non-trivial flat virtual link is the positive virtual Hopf link, $VHopf$:
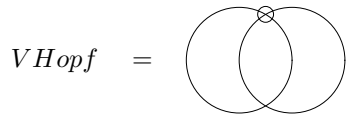
\def \VHopf{\xy 0;/r20pt/:0 *\cir(1,0){},(1,0)*\cir(1,0){},a(60)*\cir<3pt>{} \endxy}
\[
VHopf \quad=\quad \VHopf
\]

\begin{defn}
The \emph{parity} of a link is the parity of the total number of crossings between distinct
components.  The parity is an invariant of flat virtual links, because it is preserved
under the flat Reidemeister moves.
\end{defn}

For the virtual Hopf link, the parity is odd and hence $VHopf$ is non-trivial.

\begin{thm} \label{T:OCD_Equiv}
If a virtual knot is classical, then its related universe is
flat equivalent to the unknot.
\end{thm}
\begin{proof}
It is well known that every classical knot diagram can be unknotted by
choosing and switching a certain number of crossings.  The original
knot along with this unknotted diagram both share the same universe.
As a result, we can always find a sequence of flat moves taking a classical
universe to the unknot diagram.
\end{proof}

Note that a flat virtual diagram is the same as a knot universe, with the additional property
that it may also have virtual crossings.
For any oriented flat virtual diagram, there is an OCD associated with it
in the same way that a chord diagram is associated with an oriented universe.
Recall that the chords in an OCD have no arrows and are labelled with signs at their ends
corresponding to the crossing orientations of their corresponding curves.
Each flat Reidemeister move induces a related OCD move.
These OCD moves are the same as the arrow diagram moves in Figure~\ref{Fi:ADmoves} with the arrow endpoints
removed\footnote{To get an oriented chord back from an arrow, we drop the arrow tip
and place a sign on that endpoint that is opposite to the one on the basepoint.}.
This OCD approach will become useful in the discussion of filamentation invariance.

\section{Filamentations on Chord Diagrams} \label{S:Filamentations}

The notion of a filamentation\footnote{also known as bifilarations}
was first introduced in the early 1990's by Scott Carter
in~\cite{Carter1,Carter2,Carter3}
while looking at generic immersions of disks in 3 space.
In the particular case where the boundary of the disk is mapped to the boundary
of the manifold, he noted that the intersection curves would necessarily have
a total net intersection of zero.  In 2000, he suggested (private communication)
that filamentations
could be used to answer a conjecture of Kauffman's in~\cite{LouVirt} that
the following flat knot was non-trivial under flat virtual equivalence:
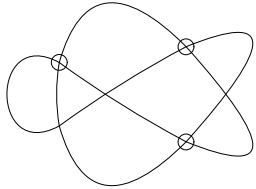
\[
{\xy /r1pc/:
,(0,-1)="b1",(0,1)="a1"*{\cir<3pt>{}}
,(4,-1.5)="a2"*{\cir<3pt>{}},(4,1.5)="b2"*{\cir<3pt>{}}
,{ \ar@{-}@/_1pt/ "a1";"a2" \ar@{-}@/^1pt/ "b1";"b2" \ar@{-}@/^1pt/ "b1";"a1" }
,"a1";"b1" **\crv{(-2.2,2.2)&(-2.2,-2.2)},"b2" **\crv{(1,4.5)}
,"a2";"b1" **\crv{(1,-4.5)},"b2" **\crv{(8.2,-3.2)},"b2" **\crv{(8.2,3.2)}
\endxy}
\]
A knot diagram with the above universe was the first example of a non-trivial
virtual knot with
trivial Jones polynomial and trivial fundamental group.

Roughly speaking, a filamentation can describe the intersection curves on an immersed
disk, much in the same way that a chord diagram describes the double points in
an immersed circle.  Consider an immersed disk which bounds a flat knot diagram.
Wherever the diagram has a flat crossing, the immersed disk will have intersection
curves emanating from the crossing.  These intersection curves will begin and end
at crossings.  A filament is a curve from the pre-image of such an intersection curve.
Thus, a filament begins at one chord endpoint and ends at another (not necessarily the same)
chord endpoint.

Our method will generalize this description so that there is no dependence on immersed
curves in space.
We describe filamentations in the purely
combinatorial context of flat virtual knots so that the existence of
a filamentation can be used to determine when a flat knot is non-trivial.

For the remainder of the paper, when we say chord diagram, we are referring generally
to OCD's and AD's.

\begin{defn}
A \emph{pairing} $\P$ on a chord diagram $\D$ is a collection of chord pairs
such that each chord in $\D$ occurs in exactly one pair in the collection.  A chord is
allowed to be paired with itself.
\end{defn}

For example, the OCD in Figure~\ref{Fi:OCDexample} has the
following possible pairings:
\begin{equation}
\begin{gathered}
\{(a,a),(b,c) \} \\
\{(b,b),(a,c) \} \\
\{(c,c),(a,b) \} \\
\{(a,a),(b,b),(c,c) \} \\
\end{gathered}
\end{equation}

\begin{defn}
A \emph{filament} $\alpha$ associated with a chord pair $(x,y)$
is a generic curve between an endpoint $X^{\epsilon}$ of $x$ and the
corresponding endpoint $Y^{-\epsilon}$ of $y$ (where $\epsilon\in\{+,-\}$).
Between endpoints, the curve must lie completely in the interior region of
one of the circles
in $\D$ and may contain a finite number of transverse self-intersections.
We orient the filament from the negative endpoint to the positive one.
\end{defn}
In general, there are an infinite number of filaments associated with a given pair.
For this paper, we will treat a filament associated with a pair as a single class,
taken up to planar isotopy of immersed curves with fixed endpoints.
\begin{defn}
The \emph{dual} of a filament $\alpha:X^-\to Y^+$, denoted $\alpha'$,
is a filament $\alpha':Y^-\to X^+$ between the two corresponding
chord endpoints $Y^-$ and $X^+$ of opposite sign.
\end{defn}
\begin{defn}
If $x$ and $y$ are distinct chords, then the two filaments associated with the pair $(x,y)$,
namely $\alpha:Y^- \to X^+$ and its dual $\alpha':X^- \to Y^+$, are called \emph{bifilaments.}
\end{defn}
\begin{defn}
A \emph{monofilament} is a filament $\alpha:X^-\to X^+$
associated with a symmetric pair $(x,x)$.  In this case, $\alpha=\alpha'$.
\end{defn}

We will refer to distinct pairs $(x,y)$ as bifilament pairs
and self-pairs $(x,x)$ as monofilament pairs.

\begin{defn}
When two distinct filaments $\alpha$ and $\beta$ intersect transversally,
they have an \emph{oriented intersection number} $\alpha\OInt\beta$.
It is calculated by looking at the local orientation (using a right-handed convention)
of the crossings between each filament: looking in the direction of $\alpha$,
if $\beta$ is directed from right to left (left to right), then
$\alpha\OInt\beta=+1$ ($=-1$).  If there are no intersections $\alpha\OInt\beta=0$.
\end{defn}
Note that the above assumes that the filaments are drawn in such a manner
that they cross exactly once, if at all.  If two filaments have more than one
transverse intersection, add up all intersection numbers between $\alpha$ and
$\beta$.  This more general description is consistent with the above definition.

Consider any pair in $\P$, and suppose that $\alpha$ is a filament
associated with that pair.
\begin{defn}
The \emph{intersection number of the filament} $\alpha$ is
\[
\IntNum{\alpha}\ =\ \sum_{\gamma\notin \{\alpha,\alpha'\} } \alpha\OInt\gamma
\]
\end{defn}
If $x$ and $y$ are distinct, and $\alpha$ and its dual $\alpha'$ are
bifilaments associated with $(x,y)$,
then the \emph{intersection number of the bifilament pair} is
\[
\IntNum{(x,y)}\ =\ \IntNum{\alpha}+\IntNum{\alpha'}\ =\
\sum_{ \gamma\notin \{\alpha,\alpha'\} } \alpha\OInt\gamma + \alpha'\OInt\gamma.
\]
Similarly, if $\alpha$ is the monofilament associated with the pair $(x,x)$,
then the \emph{intersection number of the monofilament pair} is
\[
\IntNum{(x,x)}\ =\ \IntNum{\alpha}\ =\
\sum_{ \gamma\notin \{\alpha,\alpha'\} } \alpha\OInt\gamma
\]

Note that this can be expressed as a single formula:
\begin{defn}
The \emph{intersection number} of a pair $(x,y)$, is
\[
\IntNum{(x,y)}\ =\ \sum_{\beta \in \{\alpha,\alpha'\} } \sum_{\gamma\notin \{\alpha,\alpha'\} }
\beta\OInt\gamma
\]
where $\alpha$ is a filament associated with $(x,y)$.
\end{defn}
We may need to consider multiple pairings at the same time.
If this is the case, we will use a subscript $\IntNum{\quad}_{\P}$ to specify which
pairing we are using to calculate the intersection numbers.
\begin{defn}
A \emph{filamentation} $\F$ on a Chord Diagram $\D$, is a pairing
for which the related filaments contain only transverse intersections, and
the intersection number of each pair is zero.
\end{defn}

For example, the OCD in Figure~\ref{Fi:OCDexample} has a
filamentation:
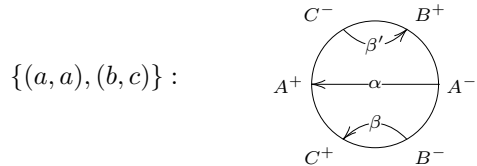
\[
\{(a,a),(b,c) \}:\qquad
\hbox{ \def\objectstyle{\scriptstyle}
\def\labelstyle{\scriptstyle}
\xy ;<24pt,0pt>:0 *\cir(1,0){}
,\ar|-{\beta'}@/_/ a(120)*+!DR{C^-};a(60)*+!DL{B^+}
,\ar|-{\beta} @/_/ a(300)*+!UL{B^-};a(240)*+!UR{C^+}
,\ar|-{\alpha} a(0)*+!L{A^-};a(180)*+!R{A^+}
\endxy}
\]

Now, we proceed with the first result.

\begin{thm} \label{T:R_fil}
If $\D$ is a chord diagram which admits a filamentation $\F$, then for any
chord diagram $\D'$ equivalent (by Reidemeister moves) to $\D$,
there is an induced filamentation $\F'$ on $\D'$.
Thus, the existence of a filamentation is an invariant of chord diagrams.
\end{thm}

\begin{proof}
Suppose $\D$ is an chord diagram and $\D'$ is equivalent to $\D$.  Then
there is a finite sequence of chord diagrams
$\{\D=\D_0,\D_1,\dots,\D_n=\D'\}$, such that
$\D_{i+1}$ and $\D_i$ differ by a single chord move.
All we need to show is that under any of the chord moves,
a filamentation can always be preserved.
This will give us a filamentation $\F_i$ on $\D_i$ induced by each move in the sequence.

In each of the following cases, we will describe how to use the
existing filamentation $\F$ to create the induced filamentation $\F'$ which
results from applying a Reidemeister move:

\begin{itemize}
\item[(R1)]
We first consider the simplest Type I cases.
%
% Case (a) and (b)
%
$${ \def\objectstyle{\scriptstyle}
\def\labelstyle{\scriptstyle}
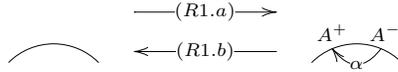
\xy ;<24pt,0pt>:
\POS(0,0) *{\xybox{*\cir(1,0){ul^dl}
    ,(.383,.924)="a" ,(-.383,.924)="b" *+!D{\phantom{A^+}}}\ \ }="d1",
\POS(4.5,0) *{\ \ \xybox{*\cir(1,0){ul^dl}
    ,(.383,.924)="a" *+!D{\ A^-} ,(-.383,.924)="b" *+!D{A^+}
    ,\ar|-{\alpha} @/^/ ,"a";"b" }}="d2",
\ar|-{(R1.a)},"d1";"d2" <20pt> \ar|-{(R1.b)},"d2";"d1" <-5pt>
\endxy}$$
If we are adding a chord, as in (R1.a), then the induced
filamentation comes from taking the old filamentation and adding
the pair $(a,a)$.  We define $\F'=\F\cup\{(a,a)\}$.  This gives us
 a monofilament $\alpha$ associated with the pair $(a,a)$.
It is clear that $\IntNum{\alpha} = \IntNum{(a,a)} = 0$ since there are no chord endpoints
on the circular arc between $A^+$ and $A^{-}$.  Thus, all
intersection numbers in the new filamentation are still zero.

If we remove a self-paired chord as in (R1.b), we again note that the
monofilament $\alpha$ associated with $(a,a)$ intersects
trivially with all other filaments associated with $\F$.
Thus we define $\F'=\F-\{(a,a)\}$, and it is clear that we still have a
filamentation.

Finally, suppose that the chord $a$ that is removed via a type I move
is paired with another chord, say $x$.

%
% Case (c)
%
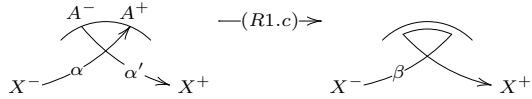
$$\hbox{ \def\objectstyle{\scriptstyle}
\def\labelstyle{\scriptstyle}
\xy ;<24pt,0pt>:
\POS(0,0)
*{\xybox{
*\cir(1,0){ul^dl}
,(.383,.924)="a" *+!D{\ A^+} ,(-.383,.924)="b" *+!D{A^-}
,(1,0)="e" *+!L{X^+} ,(-1,0)="f" *!R{X^-}
, \ar|<(.35){\alpha} @/_3pt/,"f";"a"
, \ar|<(.65){\alpha'} @/_3pt/,"b";"e"
}}="d1",
\POS(5,0)
*{\ \xybox{
*\cir(1,0){ul^dl}
,(.383,.924) *+!D{\phantom{\ A^+}} ,(-.383,.924) *+!D{\phantom{A^-}}
,(1,0)="e" *+!L{X^+} ,(-1,0)="f" *!R{X^-}
, \ar|<(.35){\beta} @{-}@/_3pt/ ,"f";(.383,.8)="a"
, \ar @{-}@/_2pt/ ,"a";(-.383,.8)="b"
, \ar @/_3pt/,"b";"e"
}}="d2",
\ar|-{(R1.c)},"d1";"d2" <10pt>
\endxy
}$$
In the case of (R1.c), the induced filamentation comes from altering the previous
pairing by replacing the pair $(x,a)$ with $(x,x)$.
That is, set $\F'=(\F-\{(x,a)\})\cup\{(x,x)\}$.  This works because we can
construct a new monofilament associated with $(x,x)$ which carries the same
intersections as the old bifilaments associated with $(x,a)$.  Consider
for example, the curve starting at $X^-$ and following the path of an
old filament $\alpha$ to $A^+$, then following along the outer circle
to the adjacent endpoint $A^-$ and finally following the path of the dual
filament $\alpha'$ to $X^+$.  If we adjust this curve slightly so that the circular
portions of the path are pushed to within the interior of the circle (as in the
above picture), we get a monofilament $\beta$ associated with $(x,x)$.  The new
curve $\beta$ contains the old intersections from both $\alpha$ and $\alpha'$.

Clearly, $\beta$ intersects the other filaments of $\F'$ only where they
coincide with the previous filaments $\alpha$ and $\alpha'$.  In addition,
any old self-intersections between $\alpha$ and $\alpha'$ will be picked up,
but we recall that intersections between dual filaments as well as
any self-intersections do not contribute to intersection numbers.  Thus, the
intersection number of the new pair is
$$\IntNum{(x,x)}_{\F'}=\IntNum{\beta}=\IntNum{\alpha}+\IntNum{\alpha'}=\IntNum{(x,a)}_{\F}=0.$$
Since all the intersections from $\F$ are preserved under the new pairing,
and no new intersections are introduced, $\F'$ is still a filamentation.
\bigskip
\item[(R2)]
We now consider the Type II moves, starting with the simplest cases.
%
% (R2.a & R2.b)
%
$$\hbox{ \def\objectstyle{\scriptstyle}
\def\labelstyle{\scriptstyle}
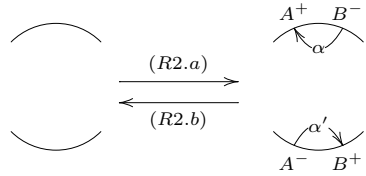
\xy 0;<24pt,0pt>:
\POS(0,0)
*{\xybox{*\cir(1,0){ul^dl} *\cir(1,0){dr^ur}
    ,(.383,.924)="a",(-.383,.924)="b" *+!D{\phantom{A^+}}
    ,(-.383,-.924)="c",(.383,-.924)="d" *+!U{\phantom{\ B^+}}} }="d1"
\POS(4,0)
*{\ \ \xybox{*\cir(1,0){ul^dl} *\cir(1,0){dr^ur}
,(.383,.924)="a" *+!D{\ B^-}
,(-.383,.924)="b" *+!D{A^+}
,(-.383,-.924)="c" *+!U{A^-}
,(.383,-.924)="d" *+!U{\ B^+}
, \ar|{\alpha} @/^8pt/ ,"a";"b", \ar|{\alpha'} @/^8pt/ ,"c";"d"
 }}="d2",
\ar^-{(R2.a)},"d1";"d2" <2pt> \ar^-{(R2.b)},"d2";"d1" <6pt>
\endxy }$$
Whenever we add chords via a Type II move as in (R2.a),
the induced filamentation comes from adding a bifilament
pair $(a,b)$ of the newly created chords.  We set
$\F'=\F\cup\{(a,b)\}$. The precise
configuration of the chord endpoints does not affect the resulting
filamentation, because in every Type II chord move, an endpoint of
one chord must be adjacent to the endpoint of opposing sign
from the other chord.  This forces the two filaments associated with
$(a,b)$ to be curves between neighboring chord endpoints along the circle.
Just as in (R1.a), each new filament will trivially intersect all
other filaments in $\F'$.  No other intersections are altered from
the old pairing, so $\F'$ is a filamentation.

Similarly, when a type II move removes a bifilament pair, as
in (R2.b), the induced filamentation is $\F'=\F-\{(a,b)\}$.

The remaining cases for type II chord removal cover the other possible
ways that the removed chords $a$ and $b$ can be paired.
$$\hbox{ \def\objectstyle{\scriptstyle}
\def\labelstyle{\scriptstyle}
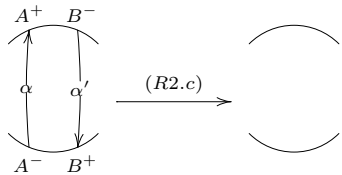
\xy 0;<24pt,0pt>:
\POS(0,0)
*{\ \ \xybox{
*\cir(1,0){ul^dl}
*\cir(1,0){dr^ur}
,(.383,.924)="a" *+!D{\ B^-}
,(-.383,.924)="b" *+!D{A^+}
,(-.383,-.924)="c" *+!U{A^-}
,(.383,-.924)="d" *+!U{\ B^+}
, \ar|{\alpha} @/^1pt/ ,"c";"b"
, \ar|{{\alpha'}} @/^1pt/ ,"a";"d"
 }}="d1",
\POS(4,0)
*{\xybox{*\cir(1,0){ul^dl}*\cir(1,0){dr^ur}
    ,(.383,.924)="a",(-.383,.924)="b" *+!D{\phantom{A^+}}
    ,(-.383,-.924)="c",(.383,-.924)="d" *+!U{\phantom{B^+}}}\ }="d2"
\ar^-{(R2.c)},"d1";"d2" <-5pt>
\endxy
}$$
Suppose $a$ and $b$ are self-paired as in (R2.c).  Note that the possible
configurations of a type II move forces the filaments to be oriented in
opposing directions, as in the picture above.  As a result, for any bifilaments
$\alpha$ and $\alpha'$ associated with (a,a) we have:
$$\IntNum{\alpha}\ = -\!\IntNum{\alpha'}$$
This also means that the net effect of $\alpha$ and $\alpha'$ on any other filament
in $\F$ cancels when computing the intersection number.
That is, for any filament $\gamma\ne\alpha,\alpha'$, we have:
\begin{eqnarray*}
\IntNum{\gamma}& = &\sum_{\delta\notin\{\gamma,\gamma'\}} \gamma\OInt\delta \\
    & = &\Bigl( \sum_{\delta\notin\{\gamma,\gamma',\alpha,\alpha'\}} \gamma\OInt\delta \Bigr) +\gamma\OInt\alpha + \gamma\OInt\alpha'\\
    & = &\sum_{\delta\notin\{\gamma,\gamma',\alpha,\alpha'\}} \gamma\OInt\delta
\end{eqnarray*}
Thus, the induced filamentation for (R2.c) is $\F'=\F-\{(a,a),(b,b)\}$.

Now, suppose that either $a$ or $b$ are paired with a third chord, $x$.
Without loss of generality we assume $a$ is self-paired and $b$ is
paired with $x$.
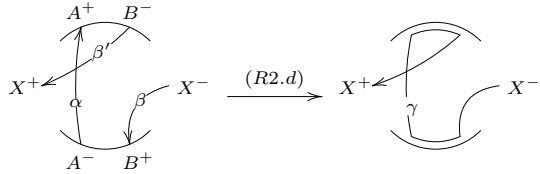
$$\hbox{ \def\objectstyle{\scriptstyle}
%
% (R2.d)
%
\xy 0;<2pc,0pc>:
\POS(0,0)
*{\xybox{ *\cir(1,0){ul^dl} *\cir(1,0){dr^ur}
,(-.383,.924)="b" *+!D{A^+} ,(-.383,-.924)="c" *+!U{A^-}
,(.383,-.924)="d" *+!U{\ \ B^+} ,(.383,.924)="a" *+!D{\ \ B^-}
,(1,0)="e" *+!L{X^-} ,(-1,0)="f" *!R{X^+}
, \ar|<(.35){\beta'} @/^2pt/ ,"a";"f"
, \ar|<(.35){\alpha} @/^2pt/,"c";"b"
, \ar|<(.40){\beta} @/_6pt/ ,"e";"d"
}\ }="d1"
\POS(5,0)
*{\ \ \xybox{ *\cir(1,0){ul^dl} *\cir(1,0){dr^ur}
,(-.383,.924)*+!D{\phantom{A^+}}, (.383,-.924) *+!U{\phantom{\ B^+}}
,(-1,0)="f" *!R{X^+} ,(1,0)="e" *+!L{X^-}
, \ar @{-}@/_6pt/ ,"e";(.383,-.8)="d"
, \ar @{-}@/^2pt/ ,"d";(-.383,-.8)="c"
, \ar|<(.30){\ \ \gamma{\phantom{'}}} @{-}@/^2pt/,"c";(-.383,.8)="b"
, \ar @{-}@/^2pt/ ,"b";(.383,.8)="a"
, \ar @/^2pt/ ,"a";"f"
 }}="d2",
\ar^-{(R2.d)},"d1";"d2" <-4pt>
\endxy
}$$
In (R2.d), the induced filamentation is
$\F'=(\F-\{(a,a),(b,x)\})\cup\{(x,x)\}$.  As before in (R1.c),
the filament associated with the new pair $(x,x)$ can be defined to be a curve
which follows the old filaments and arcs between chord endpoints so that
the locations of all intersections are preserved under the change to $\F'$.
The new monofilament $\gamma$ can start at $X^{-}$ and follow the old paths
$\beta$ to $\alpha$ to $\beta'$ (as in the picture above).
To calculate the intersection number of the new filament, we simply add the
intersection numbers of the filaments associated with the old pairs:
\begin{eqnarray*}
\IntNum{(x,x)}_{\F'}& = &\IntNum{\gamma} \\
& = &\IntNum{\alpha}+\IntNum{\beta}+\IntNum{\beta'} \\
& = &\IntNum{(a,a)}_{\F}+\IntNum{(b,x)}_{\F} \\
& = &0.
\end{eqnarray*}

Finally, suppose $a$ and $b$ are both paired with other chords.
$$
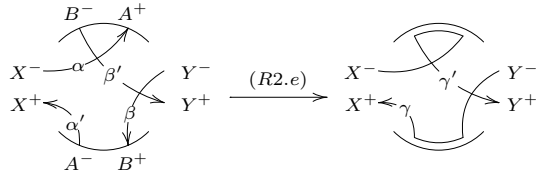
\hbox{ \def\objectstyle{\scriptstyle}
\def\labelstyle{\scriptstyle}
\xy 0;<2pc,0pc>:
\POS(0,0)
*{\xybox{ *\cir(1,0){ul^dl} *\cir(1,0){dr^ur}
,(-.383,.924)="a" *+!D{B^-} ,(.383,.924)="b" *+!D{\ A^+}
,(-.383,-.924)="c" *+!U{A^-} ,(.383,-.924)="d" *+!U{\ B^+}
,(-.95,.25)="e" *!R{X^-} ,(-.95,-.25)="f" *!R{X^+}
,(.95,.25)="g" *+!L{\ Y^-} ,(.95,-.25)="h" *+!L{\ Y^+}
, \ar|<(.35){\alpha} @/_5pt/,"e";"b"
, \ar|<(.35){\alpha'} @/_5pt/,"c";"f"
, \ar|<(.5){\beta'} @/_5pt/ ,"a";"h"
, \ar|<(.65){\beta} @/_5pt/ ,"g";"d"
}\ }="d1"
\POS(5,0)
*{\ \ \xybox{ *\cir(1,0){ul^dl} *\cir(1,0){dr^ur}
,(.383,.924) *+!D{\phantom{A^+}}
,(.383,-.924) *+!U{\phantom{B^+}}
,(-.95,.25)="e" *!R{X^-} ,(-.95,-.25)="f" *!R{X^+}
,(.95,.25)="g" *+!L{Y^-} ,(.95,-.25)="h" *+!L{Y^+}
, \ar @{-}@/_5pt/,"e";(.383,.8)="b"
, \ar @{-}@/_2pt/ ,"b";(-.383,.8)="a"
, \ar |{\gamma'}@/_5pt/ ,"a";"h"
, \ar @{-}@/_5pt/ ,"g";(.383,-.8)="d"
, \ar @{-}@/^2pt/ ,"d";(-.383,-.8)="c"
, \ar |{\ \gamma_{\ }}@/_5pt/,"c";"f"
}}="d2",
\ar^-{(R2.e)},"d1";"d2" <-4pt>
\endxy
}$$
For (R2.e), we define $\F'=(\F-\{(a,x),(b,y)\})\cup\{(x,y)\}$.
Then, we define the bifilaments associated with $(x,y)$ so that
$\gamma$ travels from $Y^-$ to $X^+$ along $\beta$ and $\alpha'$, and
the dual $\gamma'$ travels from $X^-$ to $Y^+$ along $\alpha$ and $\beta'$.
Thus,
\begin{eqnarray*}
\IntNum{(x,y)}_{\F'}& = &\IntNum{\gamma}_{\F'}+\IntNum{\gamma'}_{\F'} \\
    & = &(\IntNum{\beta}_{\F}+\IntNum{\alpha'}_{\F})+(\IntNum{\alpha}_{\F}+\IntNum{\beta'}_{\F}) \\
    & = &(\IntNum{\alpha}_{\F}+\IntNum{\alpha'}_{\F})+(\IntNum{\beta}_{\F}+\IntNum{\beta'}_{\F}) \\
    & = &\IntNum{(a,x)}_{\F}+\IntNum{(b,y)}_{\F} \\
    & = &0
\end{eqnarray*}

\item[(R3)]
Suppose we have a Type III move taking $\D$ to $\D'$.
For any possible configuration of such a move, we define the induced
filamentation to be unchanged: $\F'=\F$.

There is essentially only one $AD_3$ move to consider from
Figure~\ref{Fi:ADmoves}. This is due to the fact that the diagrams
in both versions of $AD_3$ in Figure~\ref{Fi:ADmoves} share the
same underlying OCD.  As a result, the code associated with each
move is the same (as are all of the local orientations on each
arrow endpoint)\footnote{This should not come as a surprise to the
reader since a filamentation only depends on the underlying OCD's
code, which means that it can be thought of as being associated
with the universe of the related knot diagram.}. Further, we do
not need to consider each possible arc permutation as described in
the discussion of Figure~\ref{Fi:ADmoves} in
Section~\ref{S:Definitions}. As we shall see, this is because we
can position the new filaments in such a way that they only differ
from the previous ones within a small neighborhood of each arc
involved in the move. Hence, the intersection numbers in the new
filamentation will depend only on the local filament changes near
these arcs.

% XXX -- beginning of pictures
\begin{figure}
{ \def\objectstyle{\scriptstyle}
  \def\labelstyle{\scriptstyle}
\[
\xymatrix{
    *{\xybox{0;<28pt,0pt>:0 *{\cir(1,0){}}
    \ar@{-} a(170)*+!R{A^-};a(10) *+!L{A^+\phantom{w}}
    \ar@{-} a(290)*+!UL{B^-};a(50) *+!DL{B^+}
    \ar@{-} a(130)*+!D{C^-\ \ };a(250) *+!U{C^+\ \ }
    }} \ar@{^{<}-_{>}}[rr]^-{\hbox{R3}}_{chords}
  && *{\xybox{0;<28pt,0pt>:0 *{\cir(1,0){}}
    \ar@{-} a(130)*+!DR{A^-};a(50) *+!DL{A^+}
    \ar@{-} a(250)*+!U{B^-\ \ };a(10) *+!L{B^+}
    \ar@{-} a(170)*+!R{\phantom{w}C^-};a(290) *+!UL{C^+}
    }} \\
    *{\xybox{0;<28pt,0pt>: 0 *{\cir(1,0){}} *\xycircle<15pt>{.}
        \ar @/_6pt/ a(170)*+!R{A^-};a(10) *+!L{A^+\phantom{w}}
        \ar @/^10pt/ a(290)*+!UL{B^-};a(50) *+!DL{B^+}
        \ar @/^10pt/ a(130)*+!D{C^-\ \ };a(250) *+!U{C^+\ \ }
    }} \ar@{^{<}-_{>}}[rr]^-{\hbox{R3}}_{filaments}
  && *{\xybox{0;<28pt,0pt>:0 *{\cir(1,0){}} *\xycircle<15pt>{.}
    ,a(130)*+!D{A^-\ \ };a(50)*+!DL{A^+} **\crv{(-1.0,-0.5)&(1.0,-0.5)} *\object:a(104){\dir{>}}
    ,a(250)*+!U{B^-\ \ };a(10)*+!L{B^+} **\crv{(0.7,-0.7)&(-0.4,0.9)} *\object:a(-66){\dir{>}}
    ,a(170)*+!R{\phantom{w}C^-};a(290)*+!UL{C^+} **\crv{(0.4,0.9)&(-0.7,-0.7)} *\object:a(23){\dir{>}}
        }}
    }
\]
}
\caption{\label{Fi:AD3case} An example of induced filament changes under an R3 chord move}
\end{figure}
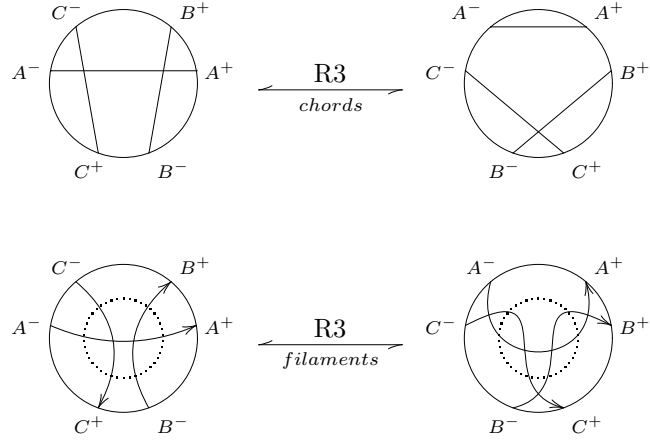
% XXX -- end of pictures
To do this, position the filaments on $\D$ so that none of the
intersections occur within a small neighborhood of the circle.  On
the new diagram $\D'$, configure the filaments as before except
within that neighborhood.  Inside this neighborhood, complete the
filament curves by crossing the two filaments which emanate from
each arc.  Note that since the $AD_3$ move switches the chord
endpoints on each arc, these new filament crossings account for
this.  We show a simple example in Figure~\ref{Fi:AD3case}. On the
left, the filaments leaving each arc don't cross until they pass
outside a neighborhood (depicted by the dotted interior circle) of
the arcs. On the right, the filaments cross as they emanate from
each arc to the dotted circle. Within the interior circle, the
filament curves on both diagrams are essentially the same.

Now, we consider the effect of adding this configuration of
three filament crossings to the original filamentation.
The general case is shown in~\ref{AD3general}.
First, suppose that the chords $a$, $b$, and $c$ are paired
with other chords, $x$, $y$, and $z$ respectively.
\begin{figure}
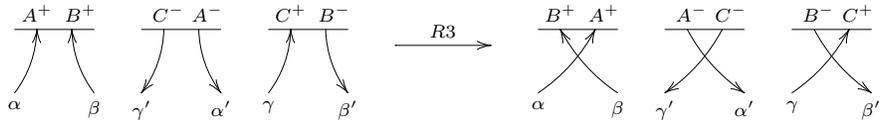

{ \def\objectstyle{\scriptstyle}
  \def\labelstyle{\scriptstyle}
\[
    \xy 0;<2pc,0pc>:0;(0.75,0)="s1",(1.25,0)="l1",(3.0,0)="s2"
        ,(0.9,0)="e1",(0.35,0)="e2",(0,1.0)="d1"
        ,"l1" **\dir{-}="a1" +"s1";p+"l1" **\dir{-}="a2" +"s1";p+"l1" **\dir{-}="a3"
        +"s2";p+"l1" **\dir{-}="b1" +"s1";p+"l1" **\dir{-}="b2" +"s1";p+"l1" **\dir{-}="b3"
        {\ar^{R3} "a3"+"s1"-(0,.25);"a3"+"s2"-"s1"-(0,.25)
        \ar @/_2pt/ "a1"-"d1"-"l1"*+!U{\alpha};"a1"-"e1" *+!D{A^+}
        \ar @/^2pt/ "a1"-"d1"*+!U{\beta};"a1"-"e2" *+!D{\ \ B^+}
        \ar @/^2pt/ "a2"-"e1" *+!D{\ C^{-}} ;"a2"-"d1"-"l1" *+!U{\gamma'}
        \ar @/_2pt/ "a2"-"e2" *+!D{\ \ A^{-}} ;"a2"-"d1" *+!U{\alpha'}
        \ar @/_2pt/ "a3"-"d1"-"l1" *+!U{\gamma};"a3"-"e1" *+!D{C^+}
        \ar @/_2pt/ "a3"-"e2" *+!D{\ \ B^{-}} ;"a3"-"d1" *+!U{\beta'}
        \ar @/_2pt/ "b1"-"d1"-"l1" *+!U{\alpha};"b1"-"e2" *+!D{\ \ A^+}
        \ar @/^2pt/ "b1"-"d1" *+!U{\beta};"b1"-"e1" *+!D{B^+}
        \ar @/^2pt/ "b2"-"e2" *+!D{\ \ C^{-}};"b2"-"d1"-"l1" *+!U{\gamma'}
        \ar @/_2pt/ "b2"-"e1" *+!D{\ A^{-}};"b2"-"d1" *+!U{\alpha'}
        \ar @/_2pt/ "b3"-"d1"-"l1" *+!U{\gamma};"b3"-"e2" *+!D{\ \ C^+}
        \ar @/_2pt/ "b3"-"e1" *+!D{\ B^{-}};"b3"-"d1" *+!U{\beta'}
        }
    \endxy
\] \bigskip
}
\caption{\label{AD3general} The general case of filament changes after an R3 chord move}
\end{figure}
The new crossings contribute as follows:
\begin{eqnarray*}
\IntNum{(a,x)}_{\F'}& = &\IntNum{(a,x)}_{\F} + \alpha\OInt\beta + \alpha'\OInt\gamma' \\
    & = &0 + 1  -1\ = \ 0\\
\IntNum{(b,y)}_{\F'}& = &\IntNum{(b,y)}_{\F} + \beta\OInt\alpha + \beta'\OInt\gamma \\
    & = &0-1+1\ = \ 0 \\
\IntNum{(c,z)}_{\F'}& = &\IntNum{(c,z)}_{\F} + \gamma\OInt\beta' + \gamma'\OInt\alpha \\
    & = &0-1+1\ = \ 0
\end{eqnarray*}
The new pairing is a filamentation whenever the old pairing is.

In fact, it does not matter how the chords are paired.  For example, suppose
we pair $a$ with $b$.  Then $x=b$, $y=a$, and from the above case, the filaments
consolidate to $\alpha=\beta'$ and $\alpha'=\beta$.
To compute $\IntNum{(a,b)}_{\F'}$, we combine the calculations of
$\IntNum{(a,x)}_{\F'}$ and $\IntNum{(b,y)}_{\F'}$ above.  First note
\[
\IntNum{(a,x)}_{\F} = \IntNum{(b,y)}_{\F}=0
\]
Although $\alpha\OInt\alpha'$ and $\alpha'\OInt\alpha$ do not contribute to
$\IntNum{(a,b)}_{\F'}$, they do sum to zero, so we will include them below to demonstrate
the similarity to $\IntNum{(a,x)}_{\F'} + \IntNum{(b,y)}_{\F'}$:
\begin{eqnarray*}
\IntNum{(a,b)}_{\F'}& = &\IntNum{(a,b)}_{\F} + \alpha\OInt\alpha' + \alpha'\OInt\gamma'
    + \alpha'\OInt\alpha + \alpha\OInt\gamma \\
    & = &0+\alpha'\OInt\gamma' + \alpha\OInt\gamma \\
    & = &-1+1 \\
    & = &0
\end{eqnarray*}
Proceed in the same fashion for all other pair choices.
\end{itemize}
This completes the proof.
\end{proof}

\begin{thm} \label{T:f_admit}
If $\D$ is a Gauss diagram which does not admit a filamentation, then the flat knot
represented by $\D$ is non-trivial.
\end{thm}

\begin{proof}
Suppose $\D$ is trivial.  Then there is a sequence of flat moves
taking $\D$ to the unknot.  Since the unknot admits a trivial filamentation,
Theorem~\ref{T:R_fil} tells us that we can find a filamentation on $\D$ by
reversing the sequence of flat moves and applying the filamentation induced
by that sequence.  Thus, if a Gauss diagram does not admit a filamentation,
it cannot represent a trivial knot.
\end{proof}

\section{An Infinite Family of Virtual Knots} \label{S:Examples}
\begin{thm} \label{T:FlatExample}
For $n\ge 2$, consider the following OCD:
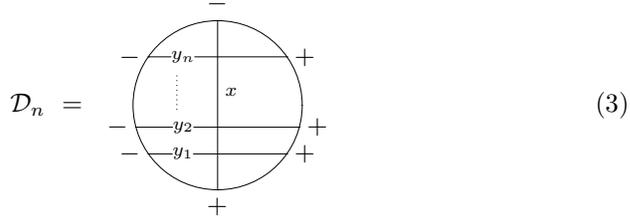
\begin{equation} \label{E:D_n}
\D_n \ = \
\xy 0;/r32pt/:0 *\cir<32pt>{}
,\ar_<(.58){x}@{-} a(270) *+!U{+};a(90)*+!D{-}
,\ar|<(.75){y_n}@{-} a(35) *+!L{+};a(145)*+!R{-}
,\ar|<(.71){y_2}@{-} a(345) *+!L{+};a(195)*+!R{-}
,\ar|<(.75){y_1}@{-} a(325) *+!L{+};a(215)*+!R{-}
,\ar@{.} (-.48,.40);(-.48,-.05)
\endxy
\end{equation}
Any flat virtual knot associated with $\D_n$ is non-trivial.
\end{thm}

\begin{proof}
Fix $n\ge 2$
Label the vertical chord $x$, and label the horizontal chords
$\{y_1,\dots,y_n\}$.
Consider the vertical chord, $x$.
We claim that for any pairing on $\D_n$, the pair including $x$
will always have a non-trivial intersection number.

First suppose $x$ is self-paired.  We see immediately that $\IntNum{(x,x)}\,=n\ne 0$,
since all of the other filaments in such a pairing must pass from left to right
across the monofilament associated with $(x,x)$.
\begin{figure}
\[
\xy 0;/r2pc/:0 *\cir<2pc>{}
,\ar@/^/ a(180) *+!R{Y_i^-};a(270)*+!U{\ X^+}
,\ar@/_/ a(90) *+!D{\ X^-};a(0)*+!L{Y_i^+}
,\ar@{.}_*\txt{$n-i$ positive\\chord endpoints}@/_2pt/ a(35)+(.3,0);a(55)+(0,.3)
,\ar@{.}_*\txt{$i-1$ negative\\chord endpoints}@/_2pt/ a(215)-(.3,0);a(235)-(0,.3)
\endxy
\]
\caption{A bifilament associated with the pair $(x,y_i)$ in $\displaystyle{D_n}$.}
\label{Fi:D_nFilaments}
\end{figure}

Suppose $x$ is paired with any of the horizontal chords, say
$y_i$. Consider Figure~\ref{Fi:D_nFilaments}.  There are exactly
$n-i$ positive chord endpoints on the arc between $Y_i^+$ and
$X^-$. Regardless of how we choose to pair the remaining chords,
these positive endpoints will have filaments ending at them.
Further, such filaments must come from the other side of the
$X^-\to Y_i^+$ filament, contributing precisely $n-1$ to the
intersection number of $(x,y_i)$. In addition, there are exactly
$i-1$ negative chord endpoints on the arc between $Y_i^-$ and
$X^+$.  Again, there must be a filament emanating from each of
these negative endpoints, and each resulting filament must cross
the $Y_i^-\to X^+$ filament.  This adds $i-1$ to the intersection
number of $(x,y_i)$.  Since we have covered all of the
intersection numbers which contribute to the bifilaments
associated with $(x,y_i)$, we get:
$$\IntNum{(x,y_i)}\,=(n-i)+(i-1)=n-1>0.$$
This means that there are no filamentations on any $\D_n$ when $n\ge 2$, so
by Theorem~\ref{T:f_admit}, any flat knot associated with $\D_n$ must be nontrivial.
\end{proof}

Theorem~\ref{T:FlatExample} gives us an infinite set of OCD's
for which any flat representative must be non-trivial.  We don't yet have a proof that
they give rise to distinct flat knots, but we believe this to be the case.
\begin{conj}
The flat knots $U_n$ in \eqref{E:D_n} are distinct for all $n\ge 1$.
\end{conj}

One interesting result we have found is a related infinite class
of virtual knot diagrams, depicted in Table~\ref{Ta:D_n}.  As we
shall see, they are all distinct.  What is even more interesting
about this class is that each knot has trivial Jones polynomial
and trivial fundamental group. This gives the first example of an
infinite class of virtual knots with trivial Jones polynomial and
trivial fundamental group.

\begin{table}
    \begin{center}
        \begin{tabular}{l||c|c|c|c|}
$n$ & $\D_n$ & $U_n$ & $\A_n$ & $K_n$ \\ \hline \hline &&&&\\
0 & {$\xy 0;<2.0pc,0pc>:0 *\cir<2pc>{},\ar@{-} a(90)*+!D{-};a(270)*+!U{+}\endxy$}
& {$\xy 0;/r1pc/:(-.2,1)="a1",(1.2,1)="b1",(0.5,.75)*{\cir<3pt>{}},(0,-.2)="a2",(1,-.2)="b2"
    ,(-.5,-.5)="a3",(1.5,-.5)="b3",(0,-.5)="c3" *{\cir<3pt>{}},(0,-1)="a4",(1,-1)="b4"
    ,\ar@/^4pt/@{-} "a1";"b2",\ar@/_4pt/@{-} "b1";"a2"
    ,\ar@{-} "a3";"b3" ,\ar@{-} "a2";"a4" ,\ar@{-} "b2";"b4"
    \hcap~{"b1"}{"b1"+(2.2,0.07)}{"b3"}{"b3"+(1.5,-.05)}
    \hcap~{"a1"}{"a1"-(2.2,-0.07)}{"a3"}{"a3"-(1.5,0.05)}|<
    \vcap~{"a4"-(0,1)}{"b4"-(0,1)}{"a4"}{"b4"} \endxy $}
& {$\xy 0;<2.0pc,0pc>:0*\cir<2pc>{},\ar a(270)*+!U{+};a(90)*+!D{\phantom{-}} \endxy$}
& {$\xy 0;/r1pc/:(-.2,1)="a1",(1.2,1)="b1",(0.5,.75)*{\cir<3pt>{}},(0,-.2)="a2",(1,-.2)="b2"
    ,(-.5,-.5)="a3",(1.5,-.5)="b3",(0,-.5)="c3" *{\cir<3pt>{}},(0,-1)="a4",(1,-1)="b4"
    ,\ar@/^4pt/@{-} "a1";"b2",\ar@/_4pt/@{-} "b1";"a2"
    ,\ar@{-} "a3";"b3" ,\ar@{-} "a2";"a4" ,\ar@{-}|!{"a3";"b3"}\hole "b2";"b4"
    \hcap~{"b1"}{"b1"+(2.2,0.07)}{"b3"}{"b3"+(1.5,-.05)}
    \hcap~{"a1"}{"a1"-(2.2,-0.07)}{"a3"}{"a3"-(1.5,0.05)}|<
    \vcap~{"a4"-(0,1)}{"b4"-(0,1)}{"a4"}{"b4"}\endxy $}
\\ &&&& \\ \hline &&&& \\
1 & ${\xy 0;<2.0pc,0pc>:0 *\cir<2pc>{},\ar@{-} a(90)*+!D{-};a(270)*+!U{+}
    ,\ar@{-} a(0)*+!L{+};a(180)*+!R{-} \endxy}$
& ${\xy 0;/r1pc/:(0,2)="a1",(1,2)="b1",(0.5,1.80)*{\cir<3pt>{}},(0,1)="a2",(1,1)="b2"
    ,(0,0)="a3",(1,0)="b3",(0,-1.3)="a4",(1,-1.3)="b4"
    ,(-.5,-1.5)="a5",(1.5,-1.5)="b5",(0,-1.5)="c5" *{\cir<3pt>{}},(0,-2)="a6",(1,-2)="b6"
    ,\ar@/^3pt/@{-} "a1";"b2",\ar@/_3pt/@{-} "b1";"a2"
    ,{\save "b2"="temp","a2";"a3":"temp"::(0,0);(1,1) **\crv{(0.5,0)&(0.5,1)},(0,1);(1,0) **\crv{(0.5,1)&(0.5,0)} \restore}
    ,{\save "b3"="temp","a3";"a4":"temp"::(0,0);(1,1) **\crv{(0.5,0)&(0.5,1)},(0,1);(1,0) **\crv{(0.5,1)&(0.5,0)},(0.5,0.5) *{\cir<3pt>{}} \restore}
    ,\ar@{-} "a5";"b5" ,\ar@{-} "a4";"a6" ,\ar@{-} "b4";"b6"
    \hcap~{"b1"}{"b1"+(3,0.05)}{"b5"}{"b5"+(2,-.05)}
    \hcap~{"a1"}{"a1"-(3,-0.05)}{"a5"}{"a5"-(2,0.05)}|<
    \vcap~{"a6"-(0,1)}{"b6"-(0,1)}{"a6"}{"b6"} \endxy}$
& ${\xy 0;<2.0pc,0pc>:0 *\cir<2pc>{},\ar a(270)*+!U{+};a(90)*+!D{\phantom{-}}
    ,\ar a(180)*+!R{-};a(0) *+!L{\phantom{+}}\endxy}$
& ${\xy 0;/r1pc/:(0,2)="a1",(1,2)="b1",(0.5,1.80)*{\cir<3pt>{}},(0,1)="a2",(1,1)="b2"
    ,(0,0)="a3",(1,0)="b3",(0,-1.2)="a4",(1,-1.2)="b4"
    ,(-.5,-1.5)="a5",(1.5,-1.5)="b5",(0,-1.5)="c5" *{\cir<3pt>{}},(0,-2)="a6",(1,-2)="b6"
    ,\ar@/^3pt/@{-} "a1";"b2",\ar@/_3pt/@{-} "b1";"a2"
    \vtwist~{"a2"}{"b2"}{"a3"}{"b3"}
    ,{\save "b3"="temp","a3";"a4":"temp"::(0,0);(1,1) **\crv{(0.5,0)&(0.5,1)},(0,1);(1,0) **\crv{(0.5,1)&(0.5,0)},(0.5,0.5) *{\cir<3pt>{}} \restore}
    ,\ar@{-} "a5";"b5" ,\ar@{-} "a4";"a6" ,\ar@{-}|!{"a5";"b5"}\hole "b4";"b6"
    \hcap~{"b1"}{"b1"+(3,0.05)}{"b5"}{"b5"+(2,-.05)}
    \hcap~{"a1"}{"a1"-(3,-0.05)}{"a5"}{"a5"-(2,0.05)}|<
    \vcap~{"a6"-(0,1)}{"b6"-(0,1)}{"a6"}{"b6"} \endxy}$
\\ &&&& \\ \hline &&&& \\
2 & ${\xy 0;<2.0pc,0pc>:0 *\cir<2pc>{}
,\ar@{-} a(270)*+!U{+};a(90)*+!D{-}
,\ar@{-} a(160)*+!R{-};a(20)*+!L{+}
,\ar@{-} a(200)*+!R{-};a(340)*+!L{+} \endxy}$
& ${\xy 0;/r1pc/:(0,3)="a1",(1,3)="b1",(0.5,2.80)*{\cir<3pt>{}},(0,2)="a2",(1,2)="b2",(0,1)="a3",(1,1)="b3"
    ,(0,0)="a4",(1,0)="b4",(0,-1)="a5",(1,-1)="b5",(0,-2.3)="a6",(1,-2.3)="b6"
    ,(-.2,-2.5)="a7",(1.2,-2.5)="b7",(0,-2.5)="c7" *{\cir<3pt>{}},(0,-3)="a8",(1,-3)="b8"
    ,\ar@/^3pt/@{-} "a1";"b2",\ar@/_3pt/@{-} "b1";"a2"
    ,{\save "b2"="temp","a2";"a3":"temp"::(0,0);(1,1) **\crv{(0.5,0)&(0.5,1)},(0,1);(1,0) **\crv{(0.5,1)&(0.5,0)} \restore}
    ,{\save "b3"="temp","a3";"a4":"temp"::(0,0);(1,1) **\crv{(0.5,0)&(0.5,1)},(0,1);(1,0) **\crv{(0.5,1)&(0.5,0)},(0.5,0.5) *{\cir<3pt>{}} \restore}
    ,{\save "b4"="temp","a4";"a5":"temp"::(0,0);(1,1) **\crv{(0.5,0)&(0.5,1)},(0,1);(1,0) **\crv{(0.5,1)&(0.5,0)} \restore}
    ,{\save "b5"="temp","a5";"a6":"temp"::(0,0);(1,1) **\crv{(0.5,0)&(0.5,1)},(0,1);(1,0) **\crv{(0.5,1)&(0.5,0)},(0.5,0.5) *{\cir<3pt>{}} \restore}
    ,\ar@{-} "a7";"b7" ,\ar@{-} "a6";"a8" ,\ar@{-} "b6";"b8"
    \hcap~{"b1"}{"b1"+(3,0.05)}{"b7"}{"b7"+(2.3,-.05)}
    \hcap~{"a1"}{"a1"-(3,-0.05)}{"a7"}{"a7"-(2.3,0.05)}|<
    \vcap~{"a8"-(0,1)}{"b8"-(0,1)}{"a8"}{"b8"} \endxy}$
& ${\xy 0;<2.0pc,0pc>:0 *\cir<2pc>{}
,\ar a(270)*+!U{+};a(90)*+!D{\phantom{-}}
,\ar a(160)*+!R{-};a(20)*+!L{\phantom{+}}
,\ar a(340)*+!L{+};a(200)*+!R{\phantom{-}} \endxy}$
& ${\xy 0;/r1pc/:(0,3)="a1",(1,3)="b1",(0.5,2.80)*{\cir<3pt>{}},(0,2)="a2",(1,2)="b2",(0,1)="a3",(1,1)="b3"
    ,(0,0)="a4",(1,0)="b4",(0,-1)="a5",(1,-1)="b5",(0,-2.0)="a6",(1,-2.0)="b6"
    ,(-.2,-2.5)="a7",(1.2,-2.5)="b7",(0,-2.5)="c7" *{\cir<3pt>{}},(0,-3)="a8",(1,-3)="b8"
    ,\ar@/^3pt/@{-} "a1";"b2",\ar@/_3pt/@{-} "b1";"a2"
    \vtwistneg~{"a2"}{"b2"}{"a3"}{"b3"}
    ,{\save "b3"="temp","a3";"a4":"temp"::(0,0);(1,1) **\crv{(0.5,0)&(0.5,1)},(0,1);(1,0) **\crv{(0.5,1)&(0.5,0)},(0.5,0.5) *{\cir<3pt>{}} \restore}
    \vtwist~{"a4"}{"b4"}{"a5"}{"b5"}
    ,{\save "b5"="temp","a5";"a6":"temp"::(0,0);(1,1) **\crv{(0.5,0)&(0.5,1)},(0,1);(1,0) **\crv{(0.5,1)&(0.5,0)},(0.5,0.5) *{\cir<3pt>{}} \restore}
    ,\ar@{-} "a7";"b7" ,\ar@{-} "a6";"a8" ,\ar@{-}|!{"a7";"b7"}\hole "b6";"b8"
    \hcap~{"b1"}{"b1"+(3,0.05)}{"b7"}{"b7"+(2.3,-.05)}
    \hcap~{"a1"}{"a1"-(3,-0.05)}{"a7"}{"a7"-(2.3,0.05)}|<
    \vcap~{"a8"-(0,1)}{"b8"-(0,1)}{"a8"}{"b8"} \endxy}$
\\ \hline &&&& \\
3 & ${\xy 0;<2.0pc,0pc>:0 *\cir<2pc>{}
,\ar@{-} a(270)*+!U{+};a(90)*+!D{-}
,\ar@{-} a(150)*+!R{-};a(30)*+!L{+}
,\ar@{-} a(0)*+!L{+};a(180)*+!R{-}
,\ar@{-} a(210)*+!R{-};a(330)*+!L{+} \endxy}$
& ${\xy 0;/r1pc/:(0,4)="a1",(1,4)="b1",(0.5,3.80)*{\cir<3pt>{}},(0,3)="a2",(1,3)="b2",(0,2)="a3",(1,2)="b3",(0,1)="a4",(1,1)="b4"
    ,(0,0)="a5",(1,0)="b5",(0,-1)="a6",(1,-1)="b6",(0,-2)="a7",(1,-2)="b7",(0,-3.0)="a8",(1,-3.0)="b8"
    ,(-.2,-3.5)="a9",(1.2,-3.5)="b9",(0,-3.5)="c9" *{\cir<3pt>{}},(0,-4)="a10",(1,-4)="b10"
    ,\ar@/^3pt/@{-} "a1";"b2",\ar@/_3pt/@{-} "b1";"a2"
    ,{\save "b2"="temp","a2";"a3":"temp"::(0,0);(1,1) **\crv{(0.5,0)&(0.5,1)},(0,1);(1,0) **\crv{(0.5,1)&(0.5,0)} \restore}
    ,{\save "b3"="temp","a3";"a4":"temp"::(0,0);(1,1) **\crv{(0.5,0)&(0.5,1)},(0,1);(1,0) **\crv{(0.5,1)&(0.5,0)},(0.5,0.5) *{\cir<3pt>{}} \restore}
    ,{\save "b4"="temp","a4";"a5":"temp"::(0,0);(1,1) **\crv{(0.5,0)&(0.5,1)},(0,1);(1,0) **\crv{(0.5,1)&(0.5,0)} \restore}
    ,{\save "b5"="temp","a5";"a6":"temp"::(0,0);(1,1) **\crv{(0.5,0)&(0.5,1)},(0,1);(1,0) **\crv{(0.5,1)&(0.5,0)},(0.5,0.5) *{\cir<3pt>{}} \restore}
    ,{\save "b6"="temp","a6";"a7":"temp"::(0,0);(1,1) **\crv{(0.5,0)&(0.5,1)},(0,1);(1,0) **\crv{(0.5,1)&(0.5,0)} \restore}
    ,{\save "b7"="temp","a7";"a8":"temp"::(0,0);(1,1) **\crv{(0.5,0)&(0.5,1)},(0,1);(1,0) **\crv{(0.5,1)&(0.5,0)},(0.5,0.5) *{\cir<3pt>{}} \restore}
    ,\ar@{-} "a9";"b9" ,\ar@{-} "a8";"a10" ,\ar@{-} "b8";"b10"
    \hcap~{"b1"}{"b1"+(3,0.05)}{"b9"}{"b9"+(2.3,-.05)}
    \hcap~{"a1"}{"a1"-(3,-0.05)}{"a9"}{"a9"-(2.3,0.05)}|<
    \vcap~{"a10"-(0,1)}{"b10"-(0,1)}{"a10"}{"b10"} \endxy}$
& ${\xy 0;<2.0pc,0pc>:0 *\cir<2pc>{}
,\ar a(270)*+!U{+};a(90)*+!D{\phantom{-}}
,\ar a(150)*+!R{-};a(30)*+!L{\phantom{+}}
,\ar a(0)*+!L{+};a(180)*+!R{\phantom{-}}
,\ar a(210)*+!R{-};a(330)*+!L{\phantom{+}} \endxy}$
& ${\xy 0;/r1pc/:(0,4)="a1",(1,4)="b1",(0.5,3.80)*{\cir<3pt>{}},(0,3)="a2",(1,3)="b2",(0,2)="a3",(1,2)="b3",(0,1)="a4",(1,1)="b4"
    ,(0,0)="a5",(1,0)="b5",(0,-1)="a6",(1,-1)="b6",(0,-2)="a7",(1,-2)="b7",(0,-3.1)="a8",(1,-3.1)="b8"
    ,(-.2,-3.5)="a9",(1.2,-3.5)="b9",(0,-3.5)="c9" *{\cir<3pt>{}},(0,-4)="a10",(1,-4)="b10"
    ,\ar@/^3pt/@{-} "a1";"b2",\ar@/_3pt/@{-} "b1";"a2"
    \vtwist~{"a2"}{"b2"}{"a3"}{"b3"}
    ,{\save "b3"="temp","a3";"a4":"temp"::(0,0);(1,1) **\crv{(0.5,0)&(0.5,1)},(0,1);(1,0) **\crv{(0.5,1)&(0.5,0)},(0.5,0.5) *{\cir<3pt>{}} \restore}
    \vtwistneg~{"a4"}{"b4"}{"a5"}{"b5"}
    ,{\save "b5"="temp","a5";"a6":"temp"::(0,0);(1,1) **\crv{(0.5,0)&(0.5,1)},(0,1);(1,0) **\crv{(0.5,1)&(0.5,0)},(0.5,0.5) *{\cir<3pt>{}} \restore}
    \vtwist~{"a6"}{"b6"}{"a7"}{"b7"}
    ,{\save "b7"="temp","a7";"a8":"temp"::(0,0);(1,1) **\crv{(0.5,0)&(0.5,1)},(0,1);(1,0) **\crv{(0.5,1)&(0.5,0)},(0.5,0.5) *{\cir<3pt>{}} \restore}
    ,\ar@{-} "a9";"b9" ,\ar@{-} "a8";"a10" ,\ar@{-}|!{"a9";"b9"}\hole "b8";"b10"
    \hcap~{"b1"}{"b1"+(3,0.05)}{"b9"}{"b9"+(2.3,-.05)}
    \hcap~{"a1"}{"a1"-(3,-0.05)}{"a9"}{"a9"-(2.3,0.05)}|<
    \vcap~{"a10"-(0,1)}{"b10"-(0,1)}{"a10"}{"b10"} \endxy}$ \\ \hline
$\vdots$ & $\vdots$ & $\vdots$ & $\vdots$ & $\vdots$
        \end{tabular}
    \end{center}
\caption{The class of virtual knots arising from $\displaystyle{D_n}$.}
\label{Ta:D_n}
\end{table}

\subsection{The Jones Polynomial via the Kauffman Bracket} \label{S:Bracket}
Many of the results stated in this section are from~\cite{LouVirt,LouVirt2}.
Another nice reference containing the known properties of the (normalized) Kauffman
bracket and the Jones polynomial is~\cite{Kawauchi}.

To compute the bracket polynomial in the case of classical knots, we start with
the skein relation:
\begin{equation} \label{E:skein}
\Bracket{\xy 0;/r1pc/:(0,1),\xoverv[2] \endxy}
\quad = \quad
A\ \Bracket{\xy 0;/r1pc/:(0,1),\huncross[2] \endxy}
\ + \
A^{-1}\ \Bracket{\xy 0;/r1pc/:(0,1) ,{\vuncross[2]}\ \endxy}
\end{equation}
Using \eqref{E:skein}, expand each crossing in a knot diagram until
we are left with a sum of collections of closed curves, called \emph{states}.
To evaluate the bracket on a state, we apply the following:
\begin{gather}
\Bracket{\xy 0;/r1pc/:0 *\cir(1,0){},a(45)*{\object:a(38){\dir{>}}}\endxy\ K\!}
\ = \ \delta\ \bracket{ K }\ = \ (-A^{2}-A^{-2})\ \bracket{ K } \\
\Bracket{\xy 0;/r1pc/:0 *\cir(1,0){},a(45)*{\object:a(37){\dir{>}}}\endxy}\ = \ 1
\end{gather}

The bracket extends naturally to the virtual category.  As above, a state of a
virtual knot diagram is the result of smoothing all real crossings.
A closed curve in such a state might still contain virtual crossings.
Ignore the virtual crossings and count the number of closed curves $||S||$ in the state.
Assign the value of $\delta^{||S||-1}$ to each state.  Since the bracket
is a regular isotopy invariant, we need to normalize it to get an invariant of virtual
knots under ambient isotopy.

\begin{defn}
Let $K$ be a virtual knot diagram.  The \emph{writhe} of $K$ is the sum
\[
w(K)=\sum_{x\in C(K)}\epsilon(x),
\]
where $C(K)$ denotes the set of all real crossings in $K$.
\end{defn}

\begin{defn}
If $K$ is a virtual knot, then the normalized bracket is given by
\[
f_{K}(A)=(-A^{3})^{-wr(K)}\bracket{K}
\]
\end{defn}

For example, we compute the bracket of the positive virtual Hopf link:
\def \VHopf#1{\xy 0;/r20pt/:0 *\cir(1,0){u^r},(1,0)*\cir(1,0){r^d},a(60)*\cir<3pt>{},#1 \endxy}
\begin{align*}
\bracket{VHopf_+}&=\Bracket{\VHopf{/u1pc/,(0,.08)*\dir{>},(1,0.08)*\dir{>}\vunder~{(1,-1)}{(0,-1)}{(1,0)}{(0,0)}}}\\
&=A\Bracket{\VHopf{\ar@{-}@/_1pt/ (0,0);(0.1,-0.4) ,\ar@{-}@/^7pt/ (0.1,-.4);(0,-1)
                 \ar@{-}@/^1pt/ (1,0);(.9,-0.4) ,\ar@{-}@/_7pt/ (.9,-.4);(1,-1)}}
+A^{-1}\Bracket{\VHopf{\vloop~{(1,-.6)}{(0,-.6)}{(1,0)}{0}, \ar@{-}@/^1pt/(0,-1);(1,-1)}} \\
&=A+A^{-1}
\end{align*}
and the normalized bracket is then
\[
f_{VHopf_+}(A)=(-A^{3})^{-(+1)}\bracket{VHopf_+}=-A^{-3}(A+A^{-1})=-A^{-2}-A^{-4}.
\]

If $K$ is a knot diagram and $K^*$ is the diagram obtained by switching all
the crossings in $K$, we have the well known property $f_{K^*}(A)=f_{K}(A^{-1})$.
This gives us the normalized bracket on the negative virtual Hopf:
\[
f_{VHopf_-}(A)=f_{VHopf_+}(A^{-1})=-A^{2}-A^{4}
\]
An immediate result is that $VHopf_+$ and $VHopf_-$ are distinct.

Another well known result is the connection between the normalized bracket
and the Jones polynomial:
\begin{thm} \label{T:JonesIsBracket}
Let $K$ be a knot and $V_K(t)$ be the Jones polynomial of $K$.  Then
\[
V_K(t)=f_K(t^{-\frac{1}{4}})
\]
\end{thm}
In the case of classical knots (of one component), the Jones polynomial is always an element of $\R[t,t^{-1}]$,
and hence the normalized bracket gives rise to polynomials in $\R[A^4,A^{-4}]$.  However,
when generalized to include virtual knots, the Jones and normalized bracket polynomials
on single component virtual knots
turn out to be in $\R[t^\frac{1}{2},t^{-\frac{1}{2}}]$ and $\R[A^2,A^{-2}]$ respectively.
This means that some virtual knots can be detected by looking for
non-integral powers of $t$ in the evaluation of the Jones polynomial (or odd powers
of $A^2$ in the normalized bracket).

\subsection{Virtual Knots Related to $\A_n$ Have Trivial Jones Polynomial}

\begin{thm}[Kauffman~\cite{LouVirt}]
The bracket is invariant under the following moves
\begin{align}
\label{E:JonesEquiv}
\Bracket{{\xy 0;/r1pc/:\htwistneg~{(0,1)}{(2,1)}{(0,-1)}{(2,-1)}\endxy} }\ &= \
\Bracket{{\xy 0;/r1pc/:(1,0) *{\cir<3pt>{}},(5,0) *{\cir<3pt>{}}
,(0,1);(2,-1) **\crv{(1,1)&(1,-1)}
,(4,1);(6,-1) **\crv{(5,1)&(5,-1)}
,(2,1);(0,-1) **\crv{(1,1)&(1,-1)}
,(6,1);(4,-1) **\crv{(5,1)&(5,-1)}
\htwistneg~{(2,1)}{(4,1)}{(2,-1)}{(4,-1)}
\endxy} } \\
\Bracket{ {\xy 0;/r1pc/:\htwist~{(0,1)}{(2,1)}{(0,-1)}{(2,-1)}\endxy} } \ &= \
\Bracket{ {\xy 0;/r1pc/:(1,0) *{\cir<3pt>{}},(5,0) *{\cir<3pt>{}}
,(0,1);(2,-1) **\crv{(1,1)&(1,-1)}
,(4,1);(6,-1) **\crv{(5,1)&(5,-1)}
,(2,1);(0,-1) **\crv{(1,1)&(1,-1)}
,(6,1);(4,-1) **\crv{(5,1)&(5,-1)}
\htwist~{(2,1)}{(4,1)}{(2,-1)}{(4,-1)}
\endxy} } \notag
\end{align}
\end{thm}

Translating the above to arrow diagrams, we have:
\begin{equation} \label{E:JonesSCD}
\Bracket{\xy 0;/r20pt/:0 *\cir(1,0){ur^ul} *\cir(1,0){dl^dr}
,(-1,0)="a"*+!DL{\epsilon},(1,0)="b"
,a(43)*{\object:a(42){\dir{>}}},a(223)*{\object:a(220){\dir{>}}}
\ar "a";"b"
\endxy \ }\ = \
\Bracket{\ \xy 0;/r20pt/:0 *\cir(1,0){ur^ul} *\cir(1,0){dl^dr}
,(-1,0)="a",(1,0)="b"*+!DR{\epsilon}
,a(43)*{\object:a(42){\dir{>}}},a(223)*{\object:a(220){\dir{>}}}
\ar "b";"a"
\endxy\ }
\end{equation}

Equation~\eqref{E:JonesSCD} implies that other than the general configuration
of the chords, the bracket depends only on the local orientations of the arrow
basepoints in an arrow diagram.  The direction of an arrow can change,
provided that the local orientations of the endpoints change with it.
Thus the relevant signs stay the same, as does the writhe.
This brings us back to what we stated in Section~\ref{S:Definitions}.
The normalized bracket (and hence the Jones polynomial) is an invariant of
signed chord diagrams.  This will be the subject of another paper.

\begin{defn}
If $V$ and $V'$ are two virtual knots for which
$f_{V}(A)=f_{V'}(A)$, then $V$ and $V'$ are said to be \emph{Jones equivalent}.
The same expression can refer to chord diagrams.
\end{defn}
There are many examples of virtual knots which are Jones equivalent.
Consider the $AD_2$ move.  Changing one arrow direction results in
\begin{equation}
\label{E:JonesEquivAD2}
\Bracket{\xy ;<20pt,0pt>:0 *\cir(1,0){ur^ul} *\cir(1,0){dl^dr}
    ,a(45) *{\object:a(42){\dir{>}}}
    ,a(225) *{\object:a(222){\dir{>}}}
    \ar_<(0.18){\epsilon} a(200);a(20)
    \ar^<(0.22){{\hbox{-}}\epsilon} a(340);a(160)
\endxy}\ =\
\Bracket{\xy ;<20pt,0pt>:0 *\cir(1,0){ur^ul} *\cir(1,0){dl^dr}
    ,a(45) *{\object:a(42){\dir{>}}}
    ,a(225) *{\object:a(222){\dir{>}}}\endxy}
\ =\
\Bracket{\xy ;<20pt,0pt>:0 *\cir(1,0){ur^ul} *\cir(1,0){dl^dr}
    ,a(45) *{\object:a(42){\dir{>}}}
    ,a(225) *{\object:a(222){\dir{>}}}
    \ar^<(0.18){\epsilon} a(170);a(10)
    \ar^<(0.22){{\hbox{-}}\epsilon} a(350);a(190)
\endxy}
\end{equation}

It is easy to construct infinite families of Jones equivalent
virtual knots. Take any arrow diagram~$\A$.  Then choose two empty
arcs within it, and apply~\eqref{E:JonesEquivAD2} an arbitrary
number of times. The $K_n$ in Table~\ref{Ta:D_n} are an example of
this.
\begin{thm} \label{T:JonesA_n}
If $K_n$ is a virtual knot associated with $\A_n$, then $f_{K_n}(A)=1$.
\end{thm}
\begin{proof}
The first two AD's, $\A_0$ and $A_1$ are representatives of the trivial knot class.
$A_0$ is a direct result of applying $AD_1$ to the unknot diagram, and
$A_1$ is the result of applying $AD_2$ to the unknot diagram.
The rest of the $\A_n$ are also Jones equivalent to the unknot, because
\begin{itemize}
\item for even $n$, the $\A_n$ are Jones equivalent to $\A_0$, and
\item for odd $n$, the $\A_n$ are Jones equivalent to $\A_1$.
\end{itemize}
This is because each $\A_n$ in Table~\ref{Ta:D_n} is Jones
equivalent to $A_{n+2}$ by a single application of
\eqref{E:JonesEquivAD2} on the horizontal arrows.  As a result,
each $K_n$ is Jones equivalent to the unknot and hence the
normalized bracket will be trivial on all of them.
\end{proof}
We should point out that Jones equivalence does not necessarily come only from
transformations of the form in~\eqref{E:JonesEquiv} and~\eqref{E:JonesSCD}.
There are also classical knots such as mutants, for example, which are Jones equivalent
to each other.
We doubt that mutants can be obtained through~\eqref{E:JonesEquiv} alone.

A well known open question is the following: are there non-trivial classical knots which
are Jones equivalent to the unknot?  In other words, does the Jones polynomial
detect knottedness for classical knots?
\begin{conj}
If $K$ is a non-trivial classical knot, then $f_{K}(A)\ne 1$.
\end{conj}

\subsection{The Fundamental Group and the Quandle} \label{S:FundGroup}
\begin{defn}
A \emph{quandle}~\cite{LouVirt} $Q$ is a non-as\-so\-ci\-a\-tive algebraic system with two binary
operations represented through the symbols~$\ur{\ }$ and~$\ul{\ }$
which satisfy the following axioms:
\begin{enumerate}
\item
For every $a\in Q$, \ $a\ul a = a$ \ and \ $a\ur a = a$.
\item
For every $a,b\in Q$, \ $a\ul b \ur b = a$ \ and \ $a\ur b \ul b = a$.
\item
For every $a,b\in Q$, there is an $x\in Q$ such that \ $x=a\ur b$ \ and \ $a=x\ul b$. \\
The left/right variant of this statement must also be true: \\
For every $a,b\in Q$, there is an $x\in Q$ such that \ $x=a\ul b$ \ and \ $a=x\ur b$.
\item
For every $a,b,c\in Q$, the following two equations hold, along with their left/right variants:
\begin{align*}
a\ur b\ur c &= a\ur c\Bur{b\ur c} \\
a\lr b\lr c &= a\Blr{c\ur b}\lr b
\end{align*}
\end{enumerate}
\end{defn}
See~\cite{LouRad} for more on the formalism of the
operator notation $\ur{\ }$ and $\ul{\ }$.  The notation is useful
because it allows us to write expressions in a non-associative setting
without requiring lots of parentheses.  The operators assume a left associative
convention, and when an expression is associated differently, the over and under
bars serve as parentheses by extending over or under the entire sub-expression
acted upon.

For example, we can think of $\ur{\ }$ as the binary operator $*$,
and $\ul{\ }$ as $\bar{*}$.  Then $a*b = a\ur b$, $a\,\bar{*}\,b=a\lr b$, and the
two basic associations are:
\begin{align*}
(a*b)*c &= a\ur b\ur c \\
a*(b*c) &= a\Bur{b\ur c} \\
\intertext{
Use the same convention for both operators.  Here are some mixed expressions:
}
(a\,\bar{*}\,b)*b = a\ul b \ur b \\
a*(b\,\bar{*}\,(c*d)) = a\Bbur{b \Bul{c\ur d}} \\
\end{align*}

The axioms of a quandle make it possible to associate quandles with knots
and links in such a way that the algebraic structure is invariant under the
Reidemeister moves.

A \emph{bridge arc} in a virtual diagram is a strand between the undercrossings
of two (possibly the same) classical crossings.
In the same sense, an \emph{arc} in a diagram is a strand
directly between any two classical crossings.
This means that a bridge arc is interpreted on an arrow diagram as an arc directly between
terminating endpoints of arrows.  An arc is a portion of the circle directly between
any two arrow endpoints.  For example, there are two arcs in $V_0$ and the
related arrow diagram $\A_0$, whereas there is only one bridge arc.

%It generalizes the fundamental group of the complement of a knot.
%The basic description is in terms of generators and equivalences on the generators.
%Given a virtual diagram, we assign generators to arcs and get equivalences from
%each classical crossing.

We define a quandle $Q(K)$ associated with a knot $K$ as follows.
Start with an oriented diagram.
Assign one generator to each bridge arc.
For each classical crossing, depending on crossing orientation,
assign an equivalence according to the following rule:
\begin{equation} \label{E:Quandle}
\xy 0;<4pc,0pc>:
\ar (0,-0.5) *+!U{a\vphantom{l}};(1,0.5) *+!D{a}
\ar|(0.5)\hole (1,-0.5) *+!U{b};(0,0.5) *+!D{c=b\ur{a}}
\ar|(0.5)\hole (3,-0.5) *+!U{a\vphantom{l}};(4,0.5) *+!D{c=a\ul{b}}
\ar (4,-0.5)*+!U{b} ;(3,0.5)*+!D{b}
\endxy
\end{equation}
We say a quandle rather than the quandle, because we refer to any quandle satisfying
these properties.
There is a universal construction~\cite{Joyce} for a quandle associated with a knot, but we
will not discuss that construction here.

Setting $b\ur{a}=bab^{-1}$ and $a\ul{b}=a^{-1}ba$ in the universal construction
gives the fundamental group of the complement of a classical knot.
We will denote $\pi_1(K)$ to be the fundamental group of a virtual knot $K$
obtained via the quandle.

\subsection{Virtual Knots Related to $\A_n$ Have Trivial Quandles}
\begin{lemma} \label{L:H_Tang}
If we evaluate the quandle on
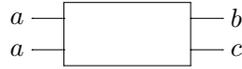
\[
\xy 0;/r1pc/:@={(0,1),(4,1),(4,-1),(0,-1)},s0="prev" @@{;"prev";**\dir{-}="prev"}
,\ar@{-} (-1,.5)*+!R{a};(0,.5)
\ar@{-} (-1,-.5)*+!R{a};(0,-.5)
\ar@{-} (5,.5)*+!L{b};(4,.5)
\ar@{-} (5,-.5)*+!L{c};(4,-.5)
\endxy
\]
where the box is replaced by a sum of the following elementary 4-tangles
$$
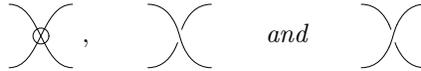
\xy 0;/r1pc/:(1,0) *{\cir<3pt>{}}
,(1,1)="a",(1,-1)="b"
,(0,1);(2,-1) **\crv{"a"&"b"}
,(2,1);(0,-1) **\crv{"a"&"b"}
\endxy
\ ,\qquad
\xy 0;/r1pc/:\htwist~{(0,1)}{(2,1)}{(0,-1)}{(2,-1)}
\endxy
\qquad\hbox{and}\qquad
\xy 0;/r1pc/:\htwistneg~{(0,1)}{(2,1)}{(0,-1)}{(2,-1)}
\endxy
$$
then $a=b=c$.  The tangle introduces no further equations.
\end{lemma}
\begin{proof}
We leave this proof as an exercise.
\end{proof}

\begin{thm} \label{T:FundGpA_n}
If\/ $K_n$ is a virtual knot associated with $\A_n$, then $\pi_1(K_n)=0$.
\end{thm}
\begin{proof}
In the right column of Table~\ref{Ta:D_n}, we have representative
virtual knots $K_n$ for each $\A_n$. By Theorem~\ref{T:AD_Equiv},
if we prove $\pi_1=0$\/ on each of these, we are done. We will
prove this via the quandle. First note that each $K_n$ in
Table~\ref{Ta:D_n} are of the general form:
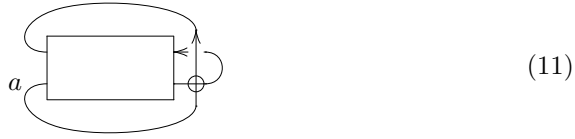
\begin{equation} \label{E:QuandleK_n}
{\xy 0;/r1pc/:@={(0,1),(4,1),(4,-1),(0,-1)},s0="prev" @@{;"prev";**\dir{-}="prev"}
,(0,.5);(4.7,1.2) **\crv{(-1,.5)&(-1,2.5)&(4.6,2)}
,(0,-.5);(4.7,-1.2) **\crv{(-1,-.5)*{a}&(-1,-2.5)&(4.6,-2)}
,(5,.5);(5,-.5) **\crv{(5.7,.5)&(5.7,-.5)}
,(4.7,-.5)*{\cir<3pt>{}}
,\ar@{-} (4,-.5);(5,-.5)
,\ar|(.3)\hole (5,.5);(4,.5)
,\ar (4.7,-1.2);(4.7,1.2)
\endxy}
\end{equation}
where the box is replaced by a horizontal sum of elementary 4-tangles.

In \eqref{E:QuandleK_n}, we have labeled the bridge arc emanating from the lower left
side of the tangle by the generator~$a$.  Note that this bridge arc both
enters and exits the left side of the tangle.
Both strands entering the left of the box are labeled with the generator~$a$.
By~\ref{L:H_Tang}, we see that the strands emanating from the right of the box
can also be labeled~$a$.
As a result, every bridge arc in the knot diagram can be labeled with
the same generator, and the quandle is trivial.
\end{proof}

\section{The Alexander Biquandle} \label{S:Biquandle}

In the construction of a quandle, we can think of each crossing as an
input/output diagram with the labels on the strands that go into the
crossing as inputs to a function, and labels on the strands that go
out of the crossing as outputs of that function.
This idea leads to the biquandle~\cite{LouFenn,LouRad,Hren1}.
With the quandle, the overcrossing strand carries a label unchanged across
the diagram, while the undercrossing strand changes its label in a manner which
depends on both input labels.

\begin{figure}
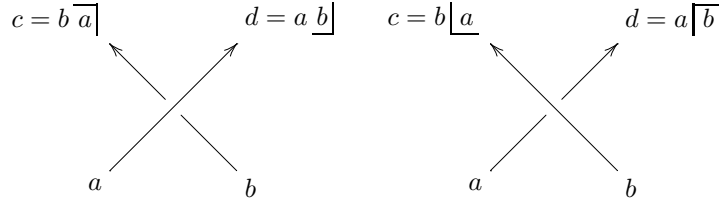

\[
\xy 0;<4pc,0pc>:
\ar ,(0,-0.5) *+!UR{a};(1,0.5) *+!DL{d=a\lr b}
\ar|(0.5)\hole ,(1,-0.5) *+!UL{b};(0,0.5) *+!DR{c=b\ur a}
\ar|(0.5)\hole ,(3,-0.5) *+!UR{a};(4,0.5) *+!DL{d=a\ul b}
\ar ,(4,-0.5)*+!UL{b} ;(3,0.5)*+!DR{c=b\ll a}
\endxy
\]
\caption{The Biquandle Operations} \label{Fi:BiquandleLabels}.
\end{figure}

In a biquandle, the overcrossing strand may also change its label.
This requires the definition of four separate functions for the
ouput strand labels, as illustrated in
 Figure~\ref{Fi:BiquandleLabels}. We indicate these functions by
the symbolism
\[
a\lr b \qquad a\ll b \qquad  a\ur b \qquad  a\ul b
\]
and view both the inputs and outputs from left to right.
Note that each of the symbols $\ur{\ }$, $\lr{\ }$, $\ll{\ }$ and $\ul{\ }$,
can be regarded as a binary operation on the underlying set of the biquandle.
Using this symbolism, the functions for the left and right crossings are
\[
R\vecfrac{a}{b}=\vecfrac{b\ur a}{a\lr b}
\qquad
L\vecfrac{a}{b}=\vecfrac{b\ll a}{a\ul b}
\]

In order for these functions to define a biquandle, they must
exhibit invariance under the Reidemeister moves.  We omit the
details here.

The Alexander biquandle is an example of a biquandle, and
we will only deal with its specific properties.

Consider any module $M$ over the ring $R=\Z[s,s^{-1},t,t^{-1}]$.
Defining the binary operations with the following equations provides us
with a biquandle structure on $M$:
\begin{alignat}{2} \label{E:ABQoperations}
    a\ur b &= ta+(1-st)b    & \qquad  a\lr b &= sa \\
    a\ul b &= \frac{1}{t}a+(1-\frac{1}{st})b    & \qquad   a\ll b &= \frac{1}{s}a
\end{alignat}
If $M$ is a free module, we call this a \emph{free Alexander biquandle}.

We associate a specific biquandle to a virtual knot diagram by taking the
free module obtained by assigning one generator for each arc and factoring
out by the submodule generated by the relations given in~\eqref{E:ABQoperations}.
We call the resulting module $ABQ(K)$ the
\emph{Alexander biquandle of the knot $K$}.

\begin{figure}
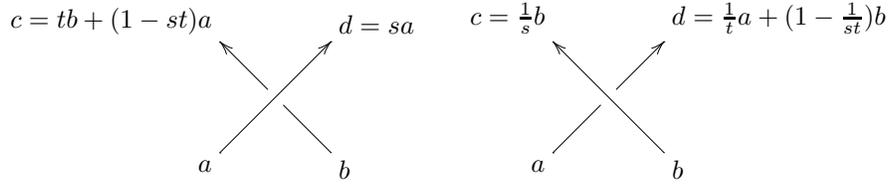

$$
\xy 0;<3.5pc,0pc>:
\ar ,(0,-0.5) *+!UR{a};(1,0.5) *+!DL{d=sa}
\ar|(0.5)\hole ,(1,-0.5) *+!UL{b};(0,0.5) *+!DR{c=tb+(1-st)a}
\ar|(0.5)\hole ,(3,-0.5) *+!UR{a};(4,0.5) *+!DL{d=\frac{1}{t}a+(1-\frac{1}{st})b}
\ar ,(4,-0.5)*+!UL{b} ;(3,0.5)*+!DR{c=\frac{1}{s}b}
\endxy
$$
% XXX -- End picture
\caption{\label{Fi:ABQrelations} The Alexander Biquandle Operations}
\end{figure}
%\begin{defn}
%The \emph{Alexander biquandle of a virtual knot $K$,} $ABQ(K)$ is a modulethe biquandle obtained by
%defining the four binary operations according to the rule in~\ref{Fi:ABQrelations}:
%\end{defn}

Note that
\begin{align}
 R\vecfrac{a}{b} &= \vecfrac{b\ur a}{a\lr b} \\
    &= \vecfrac{tb+(1-st)a}{sa} \\
    &= \begin{pmatrix} 1-st&t \\ s&0 \end{pmatrix} \vecfrac{a}{b}
\end{align}
is a linear map, and $R$ can be represented by the matrix $A$
given in Table~\ref{Ta:ABQMatrices}.  Similarly, the function $L$
can be represented by the matrix $B$ in the Table. Since we are
dealing with linear functions, we need not restrict our inputs and
outputs according to the direction of the strands.  Instead, we
can choose any two adjacent strands as inputs, and compute the
resulting function by inverting or changing the basis of the
original matrices $A$ and $B$.  To maintain a input/output
convention which is consistent with before, we will order inputs
in a counter-clockwise direction and outputs in a clockwise
direction.

Note that we have inserted a matrix $V$ in
Table~\ref{Ta:ABQMatrices}, which permutes the inputs.  This
matrix represents the virtual crossing, where labels are passed
along the strands without any changes.

\begin{table}
    \begin{align*}
        \xy 0;/r3pc/:(0,0.5) \xoverv=< \endxy: \quad A &=
            \begin{pmatrix} 1-st&t \\ s&0 \end{pmatrix} &
        \xy 0;/r3pc/:(0,0.5) \xunderv=> \endxy : \quad \widehat{A} &=
            \begin{pmatrix} 0&s \\ t&1-st \end{pmatrix} \\ \\
        \xy 0;/r3pc/:(0,0.5) \xunderv=< \endxy : \quad B &=
            \begin{pmatrix} 0&\frac{1}{s} \\ \frac{1}{t}&1-\frac{1}{st} \end{pmatrix} &
        \xy 0;/r3pc/:(0,0.5) \xoverv=> \endxy : \quad \widehat{B} &=
            \begin{pmatrix} 1-\frac{1}{st}&\frac{1}{t} \\ \frac{1}{s}&0\end{pmatrix} \\ \\
        \xy 0;/r3pc/:(0,0.5) \xoverh=> \endxy : \quad C &=
            \begin{pmatrix} 0&\frac{1}{s} \\ t&\frac{1}{s}-t \end{pmatrix} &
        \xy 0;/r3pc/:(0,0.5) \xunderh=> \endxy : \quad \widehat{C} &=
            \begin{pmatrix} \frac{1}{s}-t&t \\ \frac{1}{s}&0 \end{pmatrix} \\ \\
        \xy 0;/r3pc/:(0,0.5) \xoverh=< \endxy : \quad D &=
            \begin{pmatrix} 0&s \\ \frac{1}{t}&s-\frac{1}{t} \end{pmatrix} &
        \xy 0;/r3pc/:(0,0.5) \xunderh=< \endxy : \quad \widehat{D} &=
            \begin{pmatrix} s-\frac{1}{t}&\frac{1}{t} \\ s&0 \end{pmatrix}
    \end{align*}
    \bigskip
    \[
        \xy 0;/r1.5pc/:0*{\cir<3pt>{}},\ar@{-} (-1,-1);(1,1),\ar@{-} (-1,1);(1,-1)
        \endxy : \quad V = \begin{pmatrix} 0&1 \\ 1&0 \end{pmatrix}
    \]
    \caption{Matrices for the Alexander biquandle}
    \label{Ta:ABQMatrices}
\end{table}

The set of relations for a presentation of $ABQ(K)$
contains a generalization of the Alexander polynomial
(see~\cite{Jaeger,LouRad,Sawollek,SilverWilliams}).
\begin{defn}
The \emph{Generalized Alexander Polynomial of $K$,} $G_K(s,t)$ is the determinant
of the relation matrix from a presentation of $ABQ(K)$.  Up to multiples of $\pm s^i t^j$
for $i,j\in\Z$, it is an invariant of $K$.
\end{defn}
This polynomial is a zeroth order polynomial.
It vanishes on classical knots and links.

Recall that in the previous section, the flat knots $U_n$ are non-trivial.
The generalized Alexander polynomial provides us with the tools to
show that the related collection of knots $K_n$ are all distinct.
\begin{thm}
The virtual knots $K_n$ are all distinct for $n>0$.
\end{thm}
% \section{The Virtual Knots $K_n$ are all distinct for $n>0$.} \label{S:D_nDistinct}
\begin{proof}
We use the following model as our general form for the $K_n$:
\begin{equation} \label{E:FundGpA_n}
K_n\ =\
{\xy 0;/r1pc/:@={(1,3),(1,0),(-1,0),(-1,3)},s0="prev" @@{;"prev";**\dir{-}="prev"} *\dir{>}
,(0,1.5)*{X_n},(0,-3)*{Y}
,(1.5,-2);(1.5,-4)**\dir{-};(-1.5,-4)**\dir{-};(-1.5,-2)**\dir{-};(1.5,-2)**\dir{-} *\dir{>}
,(-.5,3);(-1.5,-3.5) **\crv{(-.5,4.5)&(-4,4.5)&(-4,-3.5)} \POS?(.74)*!DL{a}
,(-.5,0);(-1.5,-2.5) **\crv{(-.5,-1)&(-2,-.5)&(-2.7,-2.2)} \POS?(.52)*!DR{b}
,(.5,0);(1.5,-2.5) **\crv{(.5,-1)&(2,-0.5)&(2.7,-2.2)} \POS?(.52)*!DL{c}
,(.5,3);(1.5,-3.5) **\crv{(0.5,4.5)&(4,4.5)&(4,-3.5)} \POS?(.74)*!DR{d}
\endxy}
\end{equation}
The arrows show the direction we will use for the matrix insertions.
The simplified relations are:
\begin{align}
X_n\cdot\vecfrac{b}{c}&=\vecfrac{a}{d} & Y\cdot\vecfrac{b}{a}&=\vecfrac{c}{d}
\end{align}
Using the matrices from Table~\ref{Ta:ABQMatrices} to calculate
the Generalized Alexander polynomial on $K_n$, we have
\[
Y=AV=
\begin{pmatrix}
 t & 1-st \\
 0 & s
\end{pmatrix}
\]
Now set
$X=\begin{pmatrix}x_{11}&x_{12}\\ x_{21}&x_{22}\end{pmatrix}$,
so that $G_{K_n}(s,t)$ is the following determinant:
\[
G_{K_n}(s,t)\ = \
\begin{vmatrix}
1-st & t & -1 & 0 \\
s & 0 & 0 &  -1 \\
-1 & x_{11} & x_{12} & 0 \\
0 & x_{21} & x_{22} & -1
\end{vmatrix}
\]

The general form for $X_n$ depends on whether $n$ is odd or even.
\begin{equation}
X_n =
   \begin{cases}
      V(CV\widehat{C}V)^k = V(C^2)^k = VC^n,        &\text{if $n=2k$, $k\ge 0$;}    \\
      V\widehat{D}V(CV\widehat{C}V)^k = C(C^2)^k = C^n, &\text{if $n=2k+1$, $k\ge 0$.}
   \end{cases}
\end{equation}
With a little help from Maple,
\begin{equation}
C^n=
{\begin{pmatrix}
 \frac{(-t)^{n}+ts^{(1-n)}}{st+1} & \frac{-(-t)^{n}+s^{-n}}{st+1}   \\
 \frac{t(-s(-t)^{n}+s^{(1-n)})}{st+1} & \frac{ts(-t)^{n}+s^{-n}}{st+1}
\end{pmatrix}}
\end{equation}

For $n$ even
\begin{align*}
G_{K_n}(s,t)\ &= \
{\begin{vmatrix}
    1-st    & t &   -1  & 0 \\
    s   & 0 &   0   & -1    \\
    -1 & \frac{t(-s(-t)^{n}+s^{(1-n)})}{st+1} & \frac{ts(-t)^{n}+s^{-n}}{st+1} & 0  \\
     0 & \frac{(-t)^{n}+ts^{(1-n)}}{st+1} & \frac{-(-t)^{n}+s^{-n}}{st+1} & -1
\end{vmatrix}} \\ \\
 &= \frac {s^n(s^2t+1)(1-t)+s^2t^2-1+(1-s^2)t^{(1-n)}}{st+1}    \\
\intertext{ and for $n$ odd,}
G_{K_n}(s,t)\ &= \
{\begin{vmatrix}
        1-st    & t     &       -1      & 0     \\
        s       & 0     &       0       & -1    \\
        -1 & \frac{(-t)^{n}+ts^{(1-n)}}{st+1}   & \frac{-(-t)^{n}+s^{-n}}{st+1} & 0 \\
    0  & \frac{t(-s(-t)^{n}+s^{(1-n)})}{st+1} & \frac{ts(-t)^{n}+s^{-n}}{st+1}  & -1
\end{vmatrix}} \\ \\
 &= \frac{s^{(n+1)}(1-t^{2})+s^{2}t^{2}-1+(1-s^2)t^{(1-n)}}{st+1}
\end{align*}
The polynomials $s^n(s^2t+1)(1-t)$ and $s^{(n+1)}(1-t^{2})$ are distinct for any $n\ge 0$.
This proves the theorem.
\end{proof}

\section{Open Problems} \label{S:Remarks}
There are many unanswered questions in flat virtual knot theory.
In addition to the difficulty of determining when a flat knot is trivial,
it is also hard to distinguish between two flat virtual knots.
In this paper, we have given methods to determine the non-triviality of
some flat virtual knots.  We have investigated the family $U_n$
of flat virtual knots that are shadows of the knots $K_n$.
We showed in this paper that each $K_n$ is distinct from the others.
It is conjectured that the $U_n$ are all distinct as flat virtual knots.

\begin{figure}
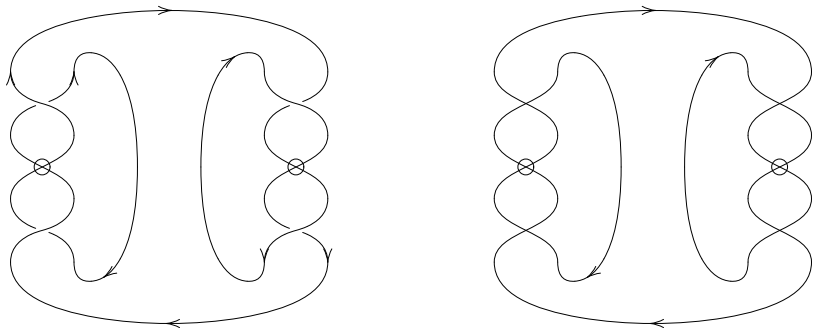

\[
{\xy 0;/r2pc/:
 (-2.5,0)="a1", (-2.5,1)="a2", (-2.5,2)="a3", (-2.5,3)="a4"
,(-1.5,0)="b1", (-1.5,1)="b2", (-1.5,2)="b3", (-1.5,3)="b4"
,(1.5,0)="c1", (1.5,1)="c2", (1.5,2)="c3", (1.5,3)="c4"
,(2.5,0)="d1", (2.5,1)="d2", (2.5,2)="d3", (2.5,3)="d4"
,{\save "b2"="temp","a2";"a3":"temp"::(0,0);(1,1) **\crv{(0.5,0)&(0.5,1)},(0,1);(1,0) **\crv{(0.5,1)&(0.5,0)},(0.5,0.5) *{\cir<3pt>{}} \restore}
,{\save "d2"="temp","c2";"c3":"temp"::(0,0);(1,1) **\crv{(0.5,0)&(0.5,1)},(0,1);(1,0) **\crv{(0.5,1)&(0.5,0)},(0.5,0.5) *{\cir<3pt>{}} \restore}
\vtwist~{"a2"}{"b2"}{"a1"}{"b1"}
\vtwist~{"c2"}{"d2"}{"c1"}{"d1"}=>
\vtwistneg~{"a4"}{"b4"}{"a3"}{"b3"}=<
\vtwistneg~{"c4"}{"d4"}{"c3"}{"d3"}
   \vcap~{"a4"+(0,2)}{"d4"+(0,2)}{"a4"}{"d4"}|<
   \vcap~{"d1"-(0,2)}{"a1"-(0,2)}{"d1"}{"a1"}|<
   \vcap~{"b4"+(0,1)}{(-0.5,4.5)}{"b4"}{(-0.5,1.5)}
   \vcap~{(-0.5,-1.5)}{"b1"-(0,1)}{(-0.5,1.5)}{"b1"}|<
   \vcap~{(0.5,4.5)}{"c4"+(0,1)}{(0.5,1.5)}{"c4"}|<
   \vcap~{"c1"-(0,1)}{(0.5,-1.5)}{"c1"}{(0.5,1.5)}
\endxy}
\qquad \qquad \qquad
{\xy 0;/r2pc/:
 (-2.5,0)="a1", (-2.5,1)="a2", (-2.5,2)="a3", (-2.5,3)="a4"
,(-1.5,0)="b1", (-1.5,1)="b2", (-1.5,2)="b3", (-1.5,3)="b4"
,(1.5,0)="c1", (1.5,1)="c2", (1.5,2)="c3", (1.5,3)="c4"
,(2.5,0)="d1", (2.5,1)="d2", (2.5,2)="d3", (2.5,3)="d4"
,{\save "b2"="temp","a2";"a3":"temp"::(0,0);(1,1) **\crv{(0.5,0)&(0.5,1)},(0,1);(1,0) **\crv{(0.5,1)&(0.5,0)},(0.5,0.5) *{\cir<3pt>{}} \restore}
,{\save "d2"="temp","c2";"c3":"temp"::(0,0);(1,1) **\crv{(0.5,0)&(0.5,1)},(0,1);(1,0) **\crv{(0.5,1)&(0.5,0)},(0.5,0.5) *{\cir<3pt>{}} \restore}
,{\save "b1"="temp","a1";"a2":"temp"::(0,0);(1,1) **\crv{(0.5,0)&(0.5,1)},(0,1);(1,0) **\crv{(0.5,1)&(0.5,0)},(0.5,0.5) \restore}
,{\save "d1"="temp","c1";"c2":"temp"::(0,0);(1,1) **\crv{(0.5,0)&(0.5,1)},(0,1);(1,0) **\crv{(0.5,1)&(0.5,0)},(0.5,0.5) \restore}
,{\save "b3"="temp","a3";"a4":"temp"::(0,0);(1,1) **\crv{(0.5,0)&(0.5,1)},(0,1);(1,0) **\crv{(0.5,1)&(0.5,0)},(0.5,0.5) \restore}
,{\save "d3"="temp","c3";"c4":"temp"::(0,0);(1,1) **\crv{(0.5,0)&(0.5,1)},(0,1);(1,0) **\crv{(0.5,1)&(0.5,0)},(0.5,0.5) \restore}
   \vcap~{"a4"+(0,2)}{"d4"+(0,2)}{"a4"}{"d4"}|<
   \vcap~{"d1"-(0,2)}{"a1"-(0,2)}{"d1"}{"a1"}|<
   \vcap~{"b4"+(0,1)}{(-0.5,4.5)}{"b4"}{(-0.5,1.5)}
   \vcap~{(-0.5,-1.5)}{"b1"-(0,1)}{(-0.5,1.5)}{"b1"}|<
   \vcap~{(0.5,4.5)}{"c4"+(0,1)}{(0.5,1.5)}{"c4"}|<
   \vcap~{"c1"-(0,1)}{(0.5,-1.5)}{"c1"}{(0.5,1.5)}
\endxy}
\]
\caption{A non-trivial virtual knot (Kishino's example)}
\label{Fi:Kishino}
\end{figure}

An intriguing example of a flat knot conjectured to be non-trivial
is the flat Kishino knot, shown in Figure~\ref{Fi:Kishino}.  This diagram is
not detected by filamentation techniques, nor any other approach that we know.
In addition, the Kishino virtual knot shown in the same figure is undetectable
by biquandles and the Jones polynomial.  It has been shown to be detected by the
3-stranded Jones polynomial~\cite{Kishino}.

\begin{figure}
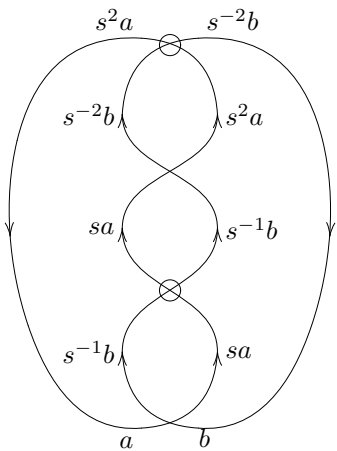

\[
{\xy 0;/r3pc/:
(0.2,3.0)="a1" *+!DR{s^2a},(0.8,3.0)="b1" *+!DL{s^{-2}b},(0.5,2.93)*{\cir<4pt>{}}
,(0,2.2)="a2" *+!R{s^{-2}b} *\dir{>},(1,2.2)="b2" *+!L{s^2a} *\dir{>}
,(0,1)="a3" *+!R{sa} *\dir{>},(1,1)="b3" *+!L{s^{-1}b} *\dir{>}
,(0,-0.3)="a4" *+!R{s^{-1}b} *\dir{>},(1,-0.3)="b4" *+!L{sa} *\dir{>}
,(0.2,-1.1)="a5" *+!UR{a},(0.8,-1.1)="b5" *!UL{b}
,\ar@/^9pt/@{-} "a1";"b2",\ar@/_9pt/@{-} "b1";"a2"
,{\save "b2"="temp","a2";"a3":"temp"::(0,0);(1,1) **\crv{(0.5,0)&(0.5,1)},(0,1);(1,0) **\crv{(0.5,1)&(0.5,0)} \restore}
,{\save "b3"="temp","a3";"a4":"temp"::(0,0);(1,1) **\crv{(0.5,0)&(0.5,1)},(0,1);(1,0) **\crv{(0.5,1)&(0.5,0)},(0.5,0.5) *{\cir<4pt>{}} \restore}
,\ar@/_9pt/@{-} "a4";"b5",\ar@/^9pt/@{-} "b4";"a5"
\hcap~{"b1"}{"b1"+(3.2,0.2)}{"b5"}{"b5"+(2.5,-.2)}|<
\hcap~{"a1"}{"a1"-(3.2,-0.2)}{"a5"}{"a5"-(2.5,0.2)}|<
\endxy}
\]
\caption{An example of how to compute the flat biquandle on a knot}
\label{Fi:FlatBQ}
\end{figure}

%Added flat BQ stuff here:
We do have a simple biquandle invariant that can detect flat links.
It is a specialization of the Alexander Biquandle that is described by the
labellings:

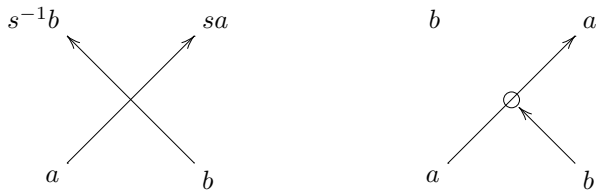
\[
\xy 0;<4pc,0pc>:
\ar ,(0,-0.5) *+!UR{a};(1,0.5) *+!DL{sa}
\ar ,(1,-0.5) *+!UL{b};(0,0.5) *+!DR{s^{-1}b}
\ar ,(3,-0.5) *+!UR{a};(4,0.5) *+!DL{a}
\ar ,(4,-0.5) *+!UL{b};(3,0.5) *+!DR{b}
,(3.5,0) *{\cir<3pt>{}}
\endxy
\]
These labellings give a module structure associated with a flat diagram
in the same way as the Alexander Biquandle.  Figure~\ref{Fi:FlatBQ} illustrates
an example link $L`$ whose flat biquandle is generated by elements $a$ and $b$,
with relations $s^2a=a$ and $s^{-2}b=b$.  Since the unlink of two
components has a module with generators $a$ and $b$ and no relations.  This
shows that $L`$ is linked (where direct parity counts fail).
Clearly, much more work remains to be accomplished in this field.

At this writing, it is not known how to extend the filamentation invariant
to links.
More generally, we would like to have more powerful combinatorial tools to
distinguish flat virtual knots and links.

\bibliographystyle{abbrv}
\bibliography{fil_paper}

\end{document}